\documentclass[final,1p,times]{elsarticle}
\usepackage[linesnumbered,ruled,boxed,vlined]{algorithm2e}
\usepackage{algorithmicx}
\usepackage{algpseudocode}
\usepackage{listings}  
\usepackage{graphicx}
\usepackage{subcaption}
\usepackage{pdflscape}

\usepackage{amssymb}
\usepackage{amsmath}
\usepackage{float}
\usepackage{amsthm}

\usepackage{multirow}
\usepackage[colorinlistoftodos]{todonotes}
\allowdisplaybreaks

\usepackage{pgf}
\usepackage{tikz}
\usetikzlibrary{arrows,automata}
\usepackage{hyperref}

\journal{}

\begin{document}

\begin{frontmatter}

\title{An adaptive large neighborhood search heuristic for the multi-port continuous berth allocation problem}

\author[label1]{Bernardo Martin-Iradi\corref{cor1}}
\address[label1]{DTU Management, Technical University of Denmark, Akademivej Building 358, 2800 Kgs. Lyngby, Denmark}

\cortext[cor1]{Corresponding author}

\ead{bmair@dtu.dk}

\author[label1]{Dario Pacino}
\ead{darpa@dtu.dk}

\author[label1]{Stefan Ropke}
\ead{ropke@dtu.dk}

\begin{abstract}
{
\footnotesize In this paper, we study a problem that integrates the vessel scheduling problem with the berth allocation into a collaborative problem denoted as the multi-port continuous berth allocation problem (MCBAP). This problem optimizes the berth allocation of a set of ships simultaneously in multiple ports while also considering the sailing speed of ships between ports. Due to the highly combinatorial character of the problem, exact methods struggle to scale to large-size instances, which points to exploring heuristic methods. We present a mixed-integer problem formulation for the MCBAP and introduce an adaptive large neighborhood search (ALNS) algorithm enhanced with a local search procedure to solve it. The computational results highlight the method's suitability for larger instances by providing high-quality solutions in short computational times. Practical insights indicate that the carriers' and terminal operators' operational costs are impacted in different ways by fuel prices, external ships at port, and the modeling of a continuous quay. 
} 
\end{abstract}


\end{frontmatter}


\section{Introduction}

The liner shipping industry is one of the major forms of international freight transportation. According to the report by \citep{unctad2020}, seaborne trade and container throughput continued growing steadily until 2019. Despite the Covid disruption during 2020, maritime trade is projected to recover and expand by 4.3 \% in 2021. 
The report also highlights that the world fleet is increasing, not only in the number of ships (more than 3 \% in 2021) but also in size. The share of the total capacity carried by mega-vessels increased from 6 \% to 40 \% in the last ten years.

This increase in demand, together with IMO's goal of reducing shipping emissions by 50 \% by 2050 \citep{IMO2018}, requires container terminals to increase capacity and improve the efficiency and sustainability of their operations.
The current growth of the vessel fleet and size directly impacts one of the most critical container terminal operations, namely the berth allocation \citep{steenken2004a}. Mathematically, this problem is denoted as the Berth Allocation Problem (BAP), which aims to assign incoming ships to berthing positions.
The BAP can assume the quay to be discrete or continuous. In the discrete version, the quay is divided into positions where each can be occupied by one ship at a time. In the continuous BAP, ships can berth at any point in the quay while respecting a safe distance from other ships.
Furthermore, the BAP can be dynamic or static. The static BAP assumes all the ships to be already at the port when the planning is done, whereas, in the dynamic version, ships can arrive at the port at different times during the planning period.
It should be noted that the dynamic BAP is still a deterministic problem. The term \textit{dynamic} refers to the different arrival times of each ship and not to the nature of the problem \citep{cordeau2005a} like in, for example, vehicle routing problems.
\begin{figure}[th]
    \centering
    \includegraphics[width=0.8\textwidth]{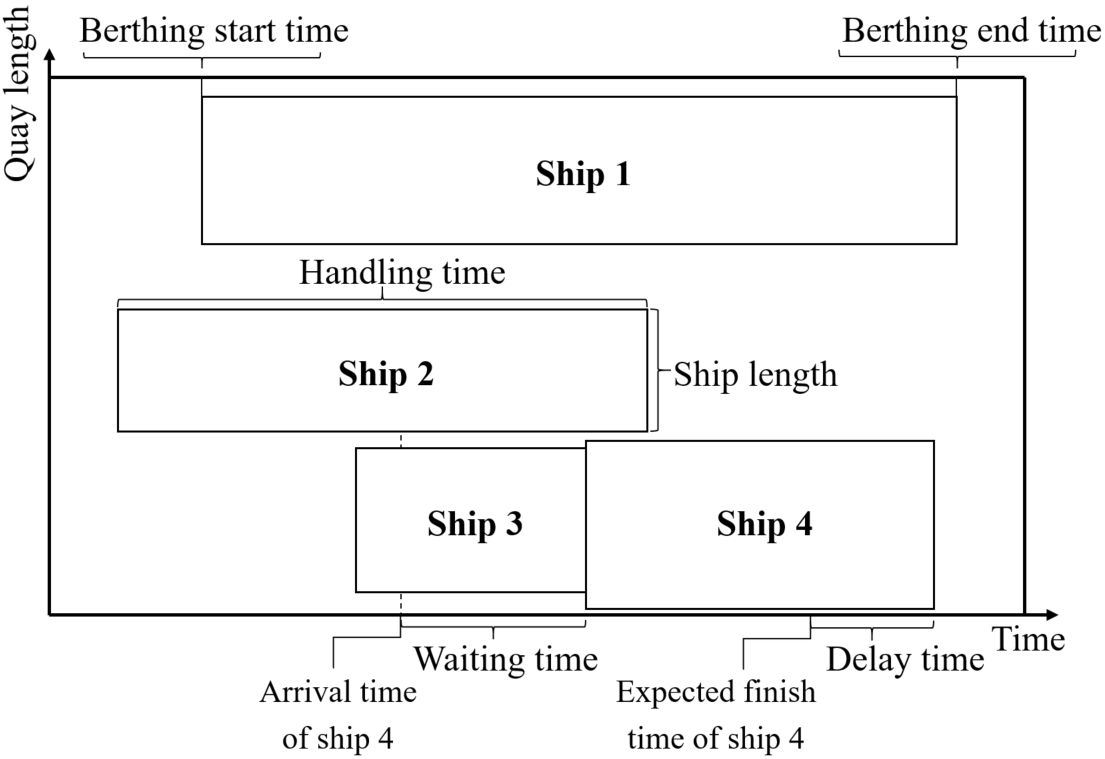}
    \caption{Example solution of the continuous and dynamic BAP for a port terminal with four vessels.}
    \label{fig:cBAPex}
\end{figure}
Figure \ref{fig:cBAPex} shows an example solution of the continuous and dynamic BAP.

Terminals optimize their berth allocation to minimize their operational costs and the time ships need to spend at the port, including waiting time, handling time, and any delays.
Due to the fierce competition between container terminals, they do not tend to share more information than is strictly required and do the planning independently from other terminals. One potential problem is that if congestion arises in a port, the affected ships can easily propagate delays to the following ports in their routes.
One way to reduce the delay is for vessels to speed up when sailing between ports. However, sailing faster results in higher fuel consumption.
This type of decision-making can be addressed by shipping line companies (i.e., carriers) in the Vessel Scheduling Problem (VSP).
The goal of the VSP is to optimize the sailing speeds between consecutive ports in the vessel's route (i.e., voyage legs). Most VSP studies aim to minimize the vessels’ fuel consumption, turnaround time at the port, and the number of vessels needed to ensure a given route frequency.
However, the VSP has its limitations. One of them is the simplistic way of modeling the berthing times of ships at port. Whereas some studies model a simplified version of berth allocation, most do not include it. Not integrating the BAP into the VSP can lead to an unrealistic or even infeasible berth allocation and, as a result, delays that ships can propagate.

A problem that integrates the berth allocation with the vessels' speed optimization was first introduced by \cite{venturini2017a} as the Multi-port Berth Allocation Problem (MBAP). This problem selects a set of ships and a set of ports that are part of their routes and simultaneously optimizes the berth allocation at all the ports, together with the sailing time between ports.
\cite{venturini2017a} studied the version of the problem with a discrete set of berthing positions.
The problem involves the joint optimization of carrier and terminal operations and relies on the a priori agreement of the vessels and ports involved. \cite{martin-iradi2022a} showed that this type of collaboration could generate cost savings for the players involved (i.e., shipping carriers and terminal operators) but also benefit the environment as fuel emissions can be reduced significantly.

As mentioned earlier in this section, the main difference between the continuous and the discrete BAP is the flexibility in the berthing positions. The set of berthing positions in the discrete BAP corresponds to a subset of those from the continuous BAP. Therefore, one can argue that modeling the quay as continuous can lead to a more resource-efficient plan, as the optimal solution of the continuous BAP is equal to or better than that of the discrete BAP. However, this potential increase in solution quality comes at the expense of higher complexity, as the solution space becomes considerably larger.

\cite{martin-iradi22b} studied the MBAP with a continuous quay (MCBAP) and highlighted the additional complexity, as the method proposed cannot scale to large instances.
This scalability issue is addressed in our study, where we employ heuristic methods that can tackle large real-world instances.

This paper makes the following four contributions:
\begin{enumerate}
    \item We define a new mixed-integer problem (MIP) formulation for the MCBAP.
    \item We present an instance generator for the MCBAP based on real-world port data, and define a set of benchmark instances that are made publically available.
    \item We implement an adaptive large neighborhood search (ALNS) method tailored to the MCBAP and enhance it with a Local Search (LS) procedure based on ejection chains.
    \item We show the viability of the ALNS method on real-size instances where it is able to find high-quality solutions faster than baseline commercial solvers.
\end{enumerate}


The remainder of this paper is structured as follows. 
Section \ref{sec:LR} comprises an extensive literature review of the MBAP together with other collaborative problems that include berth allocation or vessel scheduling. 
Section \ref{sec:prob} describes the MCBAP in detail and presents the MIP formulation.
The solution method is  described in Section \ref{sec:alns}.
Section \ref{sec:results} includes the instance generator's details and the computational study.
The conclusions and further research is summarized in Section \ref{sec:conclusions}.

\section{Literature review}\label{sec:LR}

One of the most important problems in a container terminal is the BAP, which has been studied extensively for over two decades. A survey of most of these studies is compiled in surveys by \cite{carlo2014a} and \cite{bierwirth2015a}.
\cite{lim1998a} presented one of the first formulations of the problem and showed that it is NP-hard.
Due to the additional hardness involving the BAP variant with a continuous quay, the use of heuristic methods has been predominant in the literature.
The first studies of the continuous BAP were by \cite{kim2003a} and \cite{imai2005a}, where they presented MIP formulations to the problem and solved it using heuristic and meta-heuristic algorithms such as simulated annealing.
\cite{cordeau2005a} covered both the discrete and continuous BAP and solved them using a taboo search.
\cite{guan2005a} presented both a composite heuristic and a tree search exact method and showed that both outperformed commercial solvers.
A hybrid variant between the continuous and discrete BAP was studied in \cite{kordi2016a}, where ships can only berth in a subset of positions.
Heuristic methods have been widely used when integrating the BAP with other terminal operations.
One of the main problems studied is the berth allocation and quay crane assignment problem \citep{iris2018a}.
\cite{iris2017a} present a mixed integer problem formulation with additional enhancements and implement an ALNS heuristic to solve it, whereas
\cite{cheimanoff2022a} uses a variable neighborhood search heuristic.

The VSP has also attracted significant attention in the literature.
\cite{dulebenets2019a} present a comprehensive survey about the problem and highlight the potential of collaboration and information sharing as one of the future research directions.
To the best of our knowledge, \cite{fagerholt2001a} presented the first formulation of the VSP.
Negotiating the port calls with the terminal operator \cite{dulebenets2018b} indicates that carriers and terminal operators can achieve significant savings.
A collaborative version of the VSP is presented by \cite{dulebenets2019b}, where terminal operators offer different port call durations and handling rates, leading to win-win situations.
\cite{fagerholt2010a} aim at minimizing fuel consumption by optimizing the speed in a shipping route and modeling it as a shortest path problem. The authors discretize the possible arrival times at each port to approximate the non-linear relation between fuel consumption and sailing speed. 
\cite{du2011a} and \cite{sun2018a} integrate vessel speed optimization and berth allocation by considering ships within a certain sailing distance from the port.

In the last decade, together with the increased access to data, the study of problems that require collaboration between different stakeholders (e.g., carriers and terminal operators) has become more relevant.
\cite{wang2015a} present two collaborative mechanisms that encourage sharing accurate information between carriers and terminal operators.
\cite{lalla-ruiz2016a} study the discrete BAP and present a cooperative search based on a grouping strategy where group members can only share information within the group. The collaborative berth allocation problem (CBAP) was introduced by \cite{dulebenets2018a} where a terminal planning its berth allocation can divert excessive demand to other terminals.
\cite{hellsten2020a} present an ALNS heuristic for the port scheduling problem (PSP), where the aim is to schedule feeder vessels in multi-terminal ports.
Collaboration has also been studied in disruption management. \cite{lyu2022a} present a formulation for re-planning the berth allocation and quay crane assignment and propose a heuristic method to solve it.
\cite{guo2022a} study the berth assignment and allocation problem, which integrates the BAP with the berth assignment and line clustering problem.
The first formulation of the MBAP was first introduced by \cite{venturini2017a}. It solved a dynamic and discrete BAP in multiple ports while optimizing ships' sailing speed between ports. 
\cite{martin-iradi2022a} presented a branch-and-price method for the same problem and conducted a study of the collaboration mechanism using cooperative game theory.
\cite{martin-iradi22b} extended the branch-and-price method to the MBAP with a continuous quay, the same problem of this study, and showed that exact methods are competitive for small and medium size instances but struggle to scale for larger instances. 
Recently, \cite{yu2022a} presented a genetic algorithm to solve a problem that integrates the BAP with speed optimization and vessel service differentiation to address both vertical and horizontal collaborations.

\section{Problem description}\label{sec:prob}

The MCBAP integrates operational aspects concerning terminal operators and shipping carriers. 
We consider a set of ships and a set of terminals, each of them in a different port, to optimize their operations.
Each ship visits all or a subset of the ports as a part of its route. The ships may visit the ports in different orders.
The aim of the problem is to determine the berthing position and time of the ships at each of the terminals visited.
Each terminal has a limited berthing space, given by the length of the quay. 
The service time required to load and unload the vessel is denoted as handling time and depends on the berthing position. We assume that it increases linearly with the deviation from an ideal position.
Similar to most BAP studies, the berthing time and positions of ships are subject to a set of restrictions.
Ships have a time window to be serviced also known as a port call, this is planned in advance and helps the operator to allocate berthing capacity and avoid excessive congestion. To allow for delays, the end of the time window is not strict but delays are  penalized as they require the use of unexpected resources such as more worker hours.

It is well known that the relation between sailing speed and fuel consumption is non-linear. In fact, this relation is often approximated with a cubic function as in Equation (\ref{eq:speedFuel}) \citep{venturini2017a, martin-iradi2022a}
\begin{equation}\label{eq:speedFuel}
    F(s) = (\frac{s}{s_d})^3 F_d
\end{equation}
where $s$ is the sailing speed, $s_d$ is the design speed of the ship, and $F_d$ is the fuel consumption at the design speed.
For our formulation, we discretize the set of possible sailing speeds and assume ships will sail the distance between ports at one of those speeds. Given the set of feasible sailing speeds, we can compute the corresponding set of fuel consumption rates.
This assumption ensures a linear formulation of the problem.


\begin{figure}
    \centering
    \includegraphics[width=0.84\textwidth]{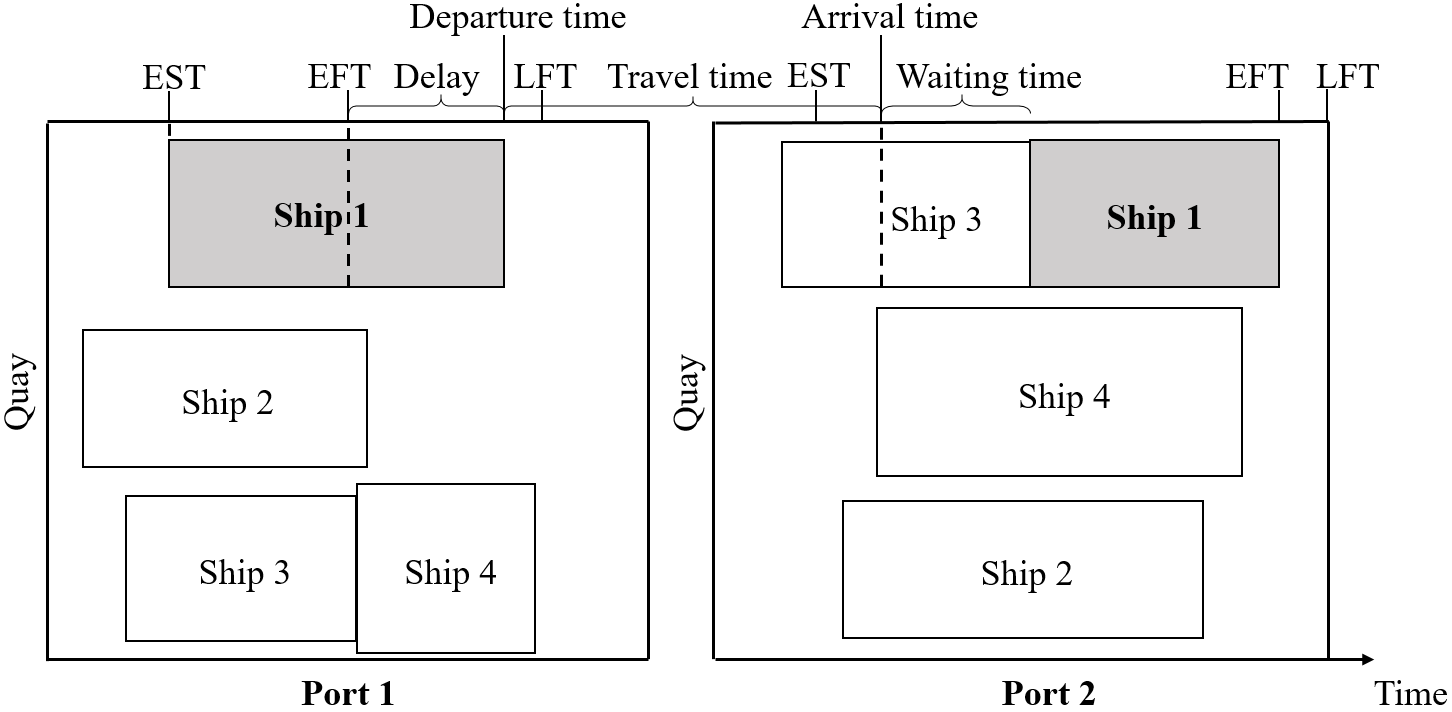}
    \caption{Example representation of a solution for the MCBAP with four ships visiting two terminals. The timeline of operations for ship 1 is defined in the top, where $EST,EFT$ and $LFT$ denote the earliest start time, the expected finish time, and the latest finish time of the ship at the port.}
    \label{fig:MCBAPex}
\end{figure}
Figure \ref{fig:MCBAPex} shows an example graphical representation of the problem, highlighting the main operational aspects of a ship (i.e., ship 1). The ship berths strictly after its earliest start time but, due to the long handling time associated with the berthing position, the service time concludes after the expected finish time. The service time after the expected finish time is computed as a delay. After the ship is serviced, it can depart towards the next port. At the time of arrival, the quay is occupied, and the ship needs to wait until ship number 3 finishes its berthing period. The time window this time is long enough to account for the waiting time, and the ship is able to finish without a delay.

\subsection{MIP formulation}

We present a new MIP formulation for the MCBAP. This formulation is based on the one for the continuous BAP from \cite{kim2003a} and the one for the discrete MBAP from \cite{venturini2017a}:
\newline

Sets and parameters:

\begin{tabular}{ p{2.0cm}  p{10.0cm} }
    $N$ & Set of all ships berthing at any of the ports. \\
    $N^* \subseteq N$ & Set of ships that we are optimizing. \\
    $\Bar{N} \subseteq N$ & Set of external ships which are considered fixed. \\
    $P$ & Set of ports. \\
    $S$ & Set of speeds. \\
    $L_p$ & Length of quay in port $p \in P$. \\
    $P_i \subseteq P$ & Set of ports planned to be visited by ship $i \in N^*$ sorted in visiting order. \\
    $C_i = \{1,...,c_i\}$ & Set of port calls for ship $i \in N$, one for each port visit. $c_i$ is the last port visit, and the value is equal to the number of port calls. \\
    $\rho_i^c$ & The port $p \in P$ corresponding to port visit $c \in C_i$ for ship $i \in N$. \\
    $N_p \subseteq N$ & Set of ships that visit port $p \in P$. \\
    $C_i^p \subseteq C_i$ & Port call positions of ship $i \in N$ visiting port $p \in P$. \\
    $x_0^{i,c}$ & The ideal berthing position for ship $i \in N^*$ at port visit $c \in C_i$ measured at the leftmost position of the ship. \\
    $EST_i^c$ & The earliest start time of berthing  for ship $i \in N^*$ at port visit $c \in C_i$. \\
    $EFT_i^c$ & The expected finish time of berthing  for ship $i \in N^*$ at port visit $c \in C_i$. \\
    $LFT_i^c$ & The latest finish time of berthing  for ship $i \in N^*$ at port visit $c \in C_i$. \\
    $\beta$ & The relative increase in handling time per unit of distance from the ideal berthing position. \\
    $\Delta^{p,p'}$ & Distance between ports $p, p' \in P$. \\
    $\Theta_s$ & Travel time per unit of distance at speed $s \in S$. \\
    $\Gamma^i_s$ & Fuel consumption per unit of distance at speed $s \in S$ for ship $i \in N^*$. \\
    $l_i$ & Length of ship $i \in N$. \\
    $F$ & Fuel cost in USD per tonne. \\
    $H$ & Cost handling time in USD per hour. \\
    $D$ & Cost of delay time in USD per hour. \\
    $I$ & Cost of waiting time in USD per hour. \\
    $U$ & Cost penalty of exceeding the latest finish time in USD. \\
\end{tabular}
\newline

Decision variables:

\begin{tabular}{ p{2.0cm}  p{10.0cm} }
    $x^{c}_{i} \in \mathbb{R}^+$ & the leftmost position of ship $i \in N$ at the quay for port visit $c \in C_i$. \\
    $y^{c}_{i} \in \mathbb{R}^+$ & the start time of berthing of ship $i \in N$ at port visit $c \in C_i$. \\
    $v_{i,s}^{c} \in \mathbb{B}$ & 1 if speed $s \in S$ is chosen by ship $i \in N^*$ to sail between port visits $c$ and $c+1$; $c \in C_i \backslash \{c_i\}$. \\
    $d_i^c \in \mathbb{R}^+$ & delay over $EFT_i^c$ for ship $i \in N^*$ at port visit $c \in C_i$. \\
    $u_i^c \in \mathbb{R}^+$ & delay over $LFT_i^c$ for ship $i \in N^*$ at port visit $c \in C_i$. \\
\end{tabular}
\newline

Auxiliary variables:

\begin{tabular}{ p{2.0cm}  p{10.0cm} }
    $\sigma^{c,c'}_{i,j} \in \mathbb{B}$ & 1 if ship $i$ is positioned left of vessel $j$ in the quay space at port visit $c \in C^p_i$ and port visit $c' \in C^p_j$ at port $p \in P$, 0 otherwise; $i,j \in N_p, i \neq j$. \\
    $\delta^{c,c'}_{i,j} \in \mathbb{B}$ & 1 if ship $i$ finishes berthing before vessel $j$ starts berthing at port visit $c \in C^p_i$ and port visit $c' \in C^p_j$ at port $p \in P$, 0 otherwise; $i,j \in N_p, i \neq j$. \\
    $r^{i,c} \in \mathbb{R}^+$ & distance between ideal and actual berthing position of ship $i \in N^*$ at port visit $c \in C_i$. \\
\end{tabular}
\newline

Dependent variables:

\begin{tabular}{ p{2.0cm}  p{10.0cm} }
    $a_i^c \in \mathbb{R}^+$ & arrival time of ship $i \in N^*$ at port visit $c \in C_i$. \\
    $h_i^c \in \mathbb{R}^+$ & handling time of ship $i \in N^*$ at port visit $c \in C_i$. \\
\end{tabular}

\begin{equation}\label{mip:obj}
    \text{min} \sum_{i \in N^*} \Big( \sum_{c \in C_i} I(y_i^c - a_i^c) + H(h_i^c) + D(d_i^c) + U(u_i^c) + \sum_{c \in C_i \backslash \{c_i\}} F(v_i^{c} \Gamma_{s}^{i} \Delta^{\rho_i^c, \rho_i^{c+1}}) \Big) 
\end{equation}

\begin{align}
    x_i^c + l_i &\leq L^p, \quad \forall i \in N_p, c \in C^p_i, p \in P \label{mip:Ql} \\
    x_{i}^c + l_{i} &\leq x^{c'}_{j}+L^p\left(1-\sigma^{c,c'}_{i,j}\right), \quad \forall p \in P, i, j \in N_p, i \neq j, c \in C^p_i, c' \in C^p_j \label{mip:leftright} \\
    y^c_{i}+h^c_{i} &\leq y^{c'}_{j}+M\left(1-\delta^{c,c'}_{i,j}\right), \quad \forall p \in P, i, j \in N_p, i \neq j, c \in C^p_i, c' \in C^p_j \label{mip:updown} \\
    \sigma^{c,c'}_{i,j}+\sigma^{c',c}_{i,j}+\delta^{c,c'}_{i,j}+\delta^{c',c}_{i,j} &\geq 1, \quad \forall i, j \in N_p, i<j, c \in C^p_i, c' \in C^p_j, c <c', p \in P \label{mip:overlap} \\
    y_i^c + h_i^c + \sum_{s \in S} v_{i,s}^{c} \Theta_s \Delta^{\rho_i^c, \rho_i^{c+1}} &= a_i^{c+1}, \quad \forall i \in N^*, c \in C_i \backslash \{c_i\} \label{mip:travelT} \\
    a_i^c &\leq y_i^c, \quad \forall i \in N^*, c \in C_i \label{mip:arrivT} \\
    EST_i^c &\leq y_i^c, \quad \forall i \in N^*, c \in C_i \label{mip:startT} \\
    y_i^c + h_i^c - EFT_i^c & \leq d_i^c \quad \forall i \in N^*, c \in C_i \label{mip:delay} \\
    y_i^c + h_i^c - LFT_i^c & \leq  u_i^c\quad \forall i \in N^*, c \in C_i \label{mip:endT} \\
     \Big(1 + \beta r^{i,c} \Big) h_0^{i,c} &= h_i^c, \quad \forall i \in N^*, c \in C_i \label{mip:hand} \\
     x_i^c - x_0^{i,c} &\leq r^{i,c}, \quad \forall i \in N^*, c \in C_i \label{mip:rplus} \\
     x_0^{i,c} - x_i^c &\leq r^{i,c}, \quad \forall i \in N^*, c \in C_i \label{mip:rminus} \\
     \sum_{s \in S} v_{i,s}^{c} &= 1, \quad \forall i \in N^*,  c \in C_i \backslash \{c_i\} \label{mip:speed} \\
     y_{i}^{c}, x_i^c & \geq 0 \quad \forall i \in N, c \in C_i\label{mip:conX} \\
      a_i^c, h_i^c, d_i^c, u_i^c, r^{i,c} & \geq 0 \quad \forall i \in N^*, c \in C_i\label{mip:conRest} \\
    v_{i,s}^{c} &\in\{0,1\} \quad \forall i \in N^*,  c \in C_i \backslash \{c_i\} \label{mip:conS} \\
    \sigma^{c,c'}_{i,j}, \delta^{c,c'}_{i,j} &\in \{0,1\}  \quad \forall i,j \in N_p, i \neq j, c \in C^p_i, c' \in C^p_j,p \in P\label{mip:conAux}
\end{align}
The set of external ships $\Bar{N}$ is considered fixed. Therefore, the corresponding set of decision variables $x_i^c, y_i^c, h_i^c, r^{i,c}$ for ships $i \in \Bar{N}$ are constant and given as input to the problem.

The objective function (\ref{mip:obj}) minimizes the operational costs of the carriers and terminal operators. This is measured as a weighted sum of the waiting time cost, handling time cost, delay cost, and fuel consumption cost.
Constraints (\ref{mip:Ql}) ensure that each ship berths within the available space. 
Constraints (\ref{mip:leftright}) and (\ref{mip:updown}) define the relative position of each pair of ships in each dimension by enabling the auxiliary variables $\sigma^{c,c'}$ and $\delta^{c,c'}$. The M value can be limited to the latest finish time of the pair of ships.
Constraints (\ref{mip:overlap}) ensure that berthing periods do not overlap in time and space.
Constraints (\ref{mip:travelT}) compute the arrival time to a port based on the sailing speed chosen to travel from the previous port.
Constraints (\ref{mip:arrivT}) and (\ref{mip:startT}) enforce that the berthing starts strictly after arrival at port and after the time window starts, respectively.
Constraints (\ref{mip:delay}) compute the delay if the expected finish time is exceeded and constraints (\ref{mip:endT}) define if the last finish time is respected.
Constraints (\ref{mip:hand}) compute the handling time for each ship and port visit while constraints (\ref{mip:rplus}) and (\ref{mip:rminus}) compute the deviation from the preferred berthing position.
Finally. constraints (\ref{mip:speed}) ensure that only one speed is chosen to sail between ports, and constraints (\ref{mip:conX}) - (\ref{mip:conAux}) define the domain of the decision variables.

\section{Solution method}\label{sec:alns}

To solve (\ref{mip:obj})-(\ref{mip:conAux}) we present an Adaptive Large Neighborhood Search (ALNS) algorithm.
The ALNS algorithm, introduced by \cite{ropke2006a}, extends the large neighborhood search method by \cite{shaw1998a}.
At each iteration, the method partially destroys and reconstructs a solution to generate a new solution. 
In our case, to destroy part of a solution, we remove the berthing time and locations of a subset of ships at a subset of ports. The combination of a scheduled berthing time and position for a ship at one of the ports in its route is denoted as a \textit{port visit}, and we will refer to this term frequently in the remainder of the paper. Additionally, in some cases, we will refer to the scheduled port visit as a \textit{rectangle}, in reference to how we can depict berthing position and time in a time-space diagram (e.g., see Figure \ref{fig:MCBAPex}).

\subsection{Construction heuristic}\label{sec:iniSol}

The ALNS requires an initial solution to start with. We present a construction heuristic process for this step that aims at finding a good initial solution.
Note that the BAP can be seen as a two-dimensional packing problem. However, in the continuous berth setting, the BAP has the increased complexity that the length of the rectangles (i.e., port visits) vary depending on the berthing location in the quay. In the case of the MCBAP, we are solving multiple  continuous BAP problems with the additional constraint that some of those berthing times depend on a sailing time. Moreover, the fact that ships follow different routes complicates the problem as greedy approaches become harder to apply.
Our construction method prioritizes reducing the delay of ships at ports. We approach this by (i) trying to place port visits early in time and close to their ideal space, therefore reducing the handling time, and (ii) by reducing "useless" space, or, in other words, placing port visits efficiently not to create empty spots in the decision space that cannot be filled by remaining port visit.
Notice that any possible solution is mathematically feasible since we allow it to exceed the latest finish time, and the time horizon is not limited. 
However, we aim to construct solutions where none of the ships exceed the $LFT$ as those can be perceived as \textit{infeasible} by the port operators and are also heavily penalized.
The method acts as a greedy heuristic, where we schedule one port visit at a time. 
The port visit to schedule is selected as the \textit{most constrained} one. 
To find it, we compute the set of feasible berthing positions and times for each ship and port visit. 
We consider a finite set of positions and times by dividing the quay into segments of a given length (e.g., 10 meters) and the planning horizon into intervals of 1 hour. 
For each time instant and segment, we compute if the ship can berth starting at that time and with its left-most side starting at the segment.
We do not count berthing times exceeding the latest finish time to measure how constrained a ship's port visit is.
From all unscheduled port visits, we define the one with the fewest possible positions as the \textit{most constrained} one.
We then schedule the port visit in one of the feasible positions. 
In fact, we do not consider the entire set but only the subset of feasible positions, where the port visit rectangle is directly adjacent to another scheduled port visit or to the limits of the decision space (i.e., the limit of the quay or planning horizon).
From this subset of positions, we select the one resulting in the minimal change to the objective function. Besides the handling and delay cost directly computed when scheduling the port visit, we need to compute fuel consumption and waiting time costs. We consider these only if the previous port visit of the ship is scheduled.

Once a port visit is scheduled, we repeat the computation and selection of the \textit{most constrained} unscheduled port visit and schedule it at the least costly \textit{efficient} position.
The procedure is described in Algorithm \ref{alg:constructionHeur}.

\begin{algorithm}
\DontPrintSemicolon
\KwData{$inst$: problem instance}
\KwResult{$sol$: a solution with all port visits scheduled}
\Begin{
$unsch \gets inst$ \tcp*{\footnotesize initialize entire set of port visits to schedule}
$sol \gets \emptyset$ \\
\While(){$unsch \neq \emptyset$}{
    \tcp*{\footnotesize sort unplanned port visits by increasing number of feasible positions}
    $unsch \gets sort(unsch)$ \\
    $toSchedule \gets popfirst(unsch)$ \tcp*{\footnotesize get first port visit from the list}
    $planned \gets false$ \\
    \tcp*{\footnotesize position at the earliest start time and closest to the ideal position}
    $pos \gets bestStartingPosition(toSchedule)$  \\
    \While(){not $planned$}{
        \If(){feasible(pos, sol)}{
            $planned \gets true$ \\
            $sol \gets plan(toSchedule, pos)$ \tcp*{\footnotesize schedule the port visit}
        }
        \Else(){
            $pos \gets updatePosition(pos,sol)$ \tcp*{\footnotesize update to next feasible position with lowest cost}
        }
    }
}
}
\caption{Construction heuristic} \label{alg:constructionHeur}
\end{algorithm}

\subsection{Removal and insertion operators}

The goal of a removal operator is to select a set of scheduled port visits to be removed from the current solution. All the operators presented in this paper select $K$ number of assignments to be removed, computed as a percentage $\rho$ of the total number of port visits to be scheduled.
It should be noted that removing port visits that are totally unrelated does not provide any potential gain. 
Therefore, a removal operator should aim at removing assignments that are related.

After applying a removal operator, the partial solution has $K$ missing port visits that need to be scheduled. They need to be assigned efficiently while respecting the other assignments 
and ensuring that the solution remains feasible. This is the goal of the insertion operators.

\subsubsection{Shaw removal}

This operator, first introduced by \cite{shaw1998a}, selects the most related pairs of assignments. To select them, we define a measure of relatedness $M_{i,j}$ between assignments $i$ and $j$ in equation \ref{eq:relate}, similar to the one presented in \cite{iris2017a}.
\begin{equation}\label{eq:relate}
    M_{i,j} = A|x_i - x_j| + B|y_i - y_j| + C|(y_i+h_i) - (y_j+h_j)|,
\end{equation}
where $x_i,y_i$ and $y_i+h_i$ are the berthing positions, berthing start time, and berthing end time of assignment $i$, respectively. $A, B$, and $C$ are custom parameters that define the importance of each of the aspects.
Observe that a lower value of $M_{i,j}$ translates into a higher level of relatedness.
To select a total of $K$ assignments, we select them following a greedy randomized criterion.
To introduce randomness in the selection of the assignments, we define a parameter $\alpha$. We sort all the port visit pairs in increasing order of $M_{i,j}$ and store them in the list $\Omega$. We then select the $i$-th element of the list applying Equation (\ref{eq:randIdx}):
\begin{equation}\label{eq:randIdx}
    i = \lceil|\Omega|\cdot p^\alpha \rceil,
\end{equation}
where $p$ is a random number $[0,1)$. Note that if $\alpha=1$, the selection is completely random, but as the value of $\alpha$ increases, the resulting value has a more deterministic behavior. The element selected will consist of two port visits to be removed. The selection process continues until $K$ port visits are removed.
Note that this method differs from the original method from \cite{shaw1998a} in that the subsequent pairs do not necessarily need to be related with the first pair selected.

\subsubsection{Time and space-relatedness removal}
This removal uses a different relatedness measure. 
We first sort all port visits by cost. 
The cost $B_i^c$ of port visit $c \in C_i$ for ship $i \in N^*$ is defined in Equation (\ref{eq:costPV}). It is measured by the ship's waiting, handling, and delay time at the port visit, plus half of the fueling costs from sailing from the previous port (if any) and to the next port (if any).
\begin{equation}\label{eq:costPV}
    B_i^c = H h_i^c + D d_i^c + I (y_i^c - a_i^c) + \frac{F^c_i}{2}
\end{equation}
$F^c_i$ is the fuel costs associated with the previous and next port visits if any. For example, if the ship sails from a previous port visit $c_p$ to port visit $c$, and then continues to the next port visit $c_n$, then the fuel costs are computed as in Equation (\ref{eq:fuelCosts}).
\begin{equation}\label{eq:fuelCosts}
    F^c_i = F(v_i^{c_p,c} \Gamma_{s}^{i} \Delta^{\rho(c_p), \rho(c)}) + F(v_i^{c,c_n} \Gamma_{s}^{i} \Delta^{\rho(c), \rho(c_n)})
\end{equation}
In the case that port visit $c$ is the first or last port visit in the route for the ship, the corresponding missing sailing leg is removed from the fuel cost computation.

We then select the $i$\textsuperscript{th} most expensive assignment applying Equation (\ref{eq:randIdx}) and remove all \textit{neighbor} assignments. We define as \textit{neighbors} all the assignments that are within a \textit{distance} of the assignment. We consider the \textit{distance} in both time and space. If an assignment is depicted as a rectangle in a time-space diagram of the port, the \textit{neighbor} area represents the one that overlaps in time or space with it. All other assignments that overlap partially or completely with the neighbor area are considered neighbors and removed. We then select the most expensive assignment and remove all neighbor assignments. We repeat the process until $K$ assignments are removed.

\begin{figure}
     \centering
     \begin{subfigure}[b]{0.45\textwidth}
         \centering
         \includegraphics[width=\textwidth]{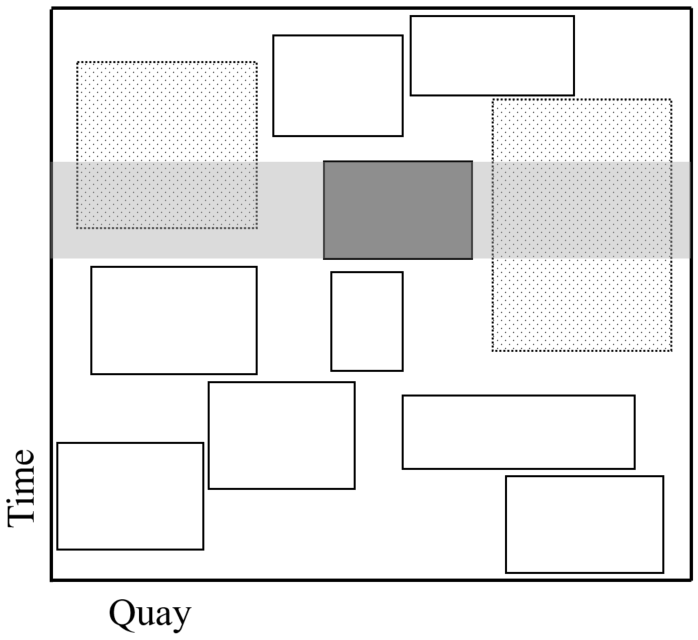}
         \caption{Time-related port visits}
         \label{fig:y equals x}
     \end{subfigure}
     \hfill
     \begin{subfigure}[b]{0.45\textwidth}
         \centering
         \includegraphics[width=\textwidth]{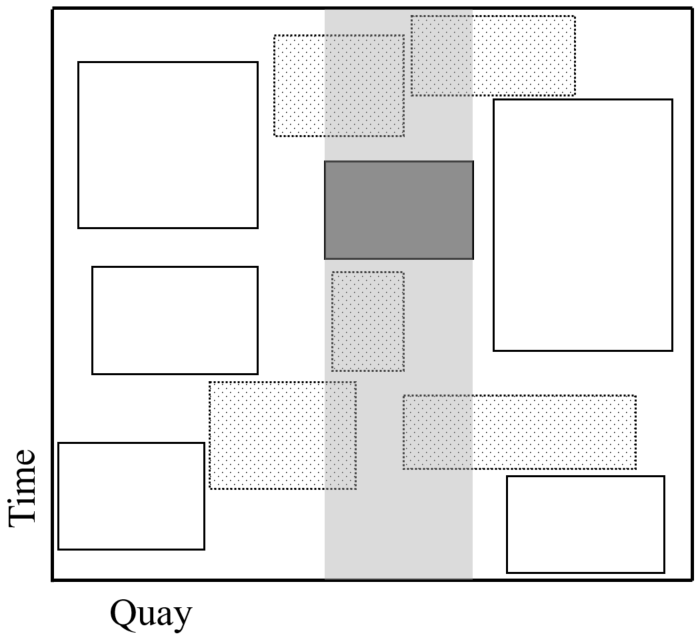}
         \caption{Space-related port visits}
         \label{fig:three sin x}
     \end{subfigure}
        \caption{Neighbor port visits in time and space for a given port visit in dark grey.}
        \label{fig:timeSpaceRelated}
\end{figure}
Figure \ref{fig:timeSpaceRelated} shows an example of neighbor port visits in time and space.
Depending on the dimension considered we define the two removal operators as \textit{cost-time removal} and \textit{cost-space removal}.

\subsubsection{Random removal}
We also consider a fully randomized destroy operator. It randomly selects $K$ assignments to be removed. The goal of this operator is not  to select relevant port visits to remove but rather to help diversify the search.

\subsubsection{Randomized greedy insertion}\label{sec:greedyRep}

This method follows the same procedure as the construction heuristic with the addition of a randomized component when selecting the port visit to schedule at each step.

All unplanned port visits are sorted based on the number of available insertion positions.
An available insertion position is one that maintains a feasible solution. For instance, the port visit needs to ensure that the previous, or following port visits, are connected through a feasible sailing speed if any of these are already scheduled.
We select the port visit using a randomization parameter $\gamma$ in the same way that $\alpha$ is used in Equation \ref{eq:randIdx}. This prioritizes the port visits with fewer available insertion positions.
The selected port visit
is scheduled in the position that increases the objective function the least (i.e., lowest cost).
The process iterates by recalculating the new number of insertion positions for the remaining port visits.

\subsubsection{$\kappa$-regret insertion}

This insertion method is based on the \textit{regret-k} heuristic presented in \cite{potvin1993a}. This method has an additional \textit{look-ahead} component compared to a basic greedy heuristic.
For each of the port visits, we compute the $\kappa$ best scheduling positions, and we then measure the \textit{regret} cost for each of them as the difference between the best and $\kappa$-best positions. The one with the highest regret cost becomes the next port visit to plan.
The process is described in Algorithm \ref{alg:kRegretRep}.

\begin{algorithm}
\DontPrintSemicolon
\KwData{$sol, unsch, \kappa$: partially destroyed solution, set of port visits to schedule, and the parameter $\kappa$}
\KwResult{$sol$: repaired solution with all port visits scheduled.}
\Begin{

\While(){$unsch \neq \emptyset$}{
    $order \gets \emptyset$ \tcp*{\small initialize empty list}
    \For(){$portVisit \in unsch$}{
        $[pos] \gets findBestPositions(\kappa)$ \tcp*{\small compute $\kappa$ best insert positions}
        $regretCost \gets c(pos[\kappa] - pos[1])$ \tcp*{\small compute regret cost}
        $order \gets sortList(portVisit, regretCost)$ \tcp*{\small update list by regret cost}
    }
    $sol \gets plan(order[1])$ \tcp*{\small plan selected port visit}
    $unsch \gets pop(order[1])$ \tcp*{\small update set of unplanned port visits}
}
}
\caption{$\kappa$-regret insertion}\label{alg:kRegretRep}
\end{algorithm}

\subsubsection{Packing greedy insertion}

This insertion method is similar to the randomized greedy insertion described in Section \ref{sec:greedyRep}. The main difference is the position where the port visits are planned.
Scheduling the port visits in a position with lower objective value can lead to the creation of empty spaces and, therefore, to inefficient use of the decision space. 
This method restricts the set of possible insertion positions to the ones \textit{strictly adjacent} to other scheduled ships, or to the limits of the quay or planning horizon.
\begin{figure}
    \centering
    \includegraphics[width=0.4\textwidth]{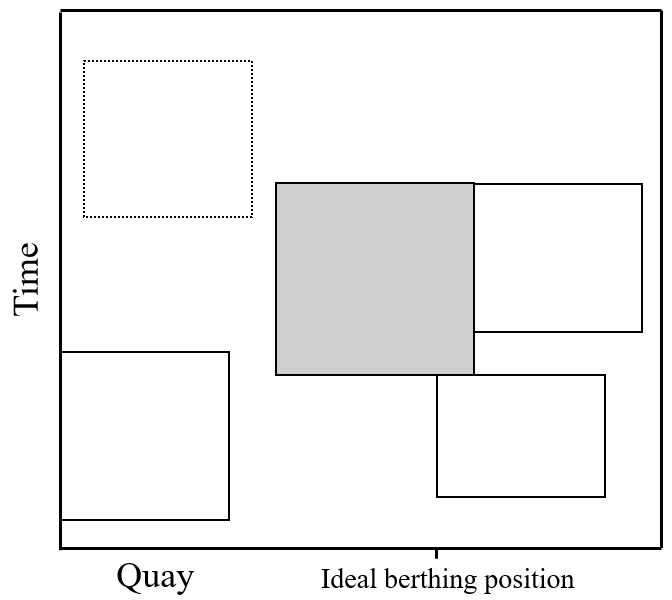}
    \caption{Graphical representation of example positions (continuous line) \textit{strictly adjacent} to the grey ship or the quay space. The position represented with a dashed line is not part of the set of positions as it is not adjacent to another planned ship or the boundaries of the decision space.}
    \label{fig:EfficPackEx}
\end{figure}
By \textit{strictly adjacent}, we mean that the port visit to schedule needs to berth strictly next to another ship during at least one interval of time (e.g., one hour) or berth strictly before (or after) another ship with at least one quay segment in common. Also, we consider berthing positions where one of the sides is at one end of the quay, or if the berthing period starts or ends at the earliest and latest possible berthing time, respectively. Figure \ref{fig:EfficPackEx} shows some example positions considered.

\subsubsection{Arrival greedy insertion}

This method is identical to the one presented in Section \ref{sec:greedyRep} with the only difference that instead of sorting the unplanned port visits by increasing the number of feasible insertion positions, we sort the unplanned port visits by the earliest possible arrival time.
One of the main goals of this method is to schedule port visits earlier, at the expense of a potentially higher cost, in order to increase the number of possible insertion positions for the remaining unplanned port visits.

\subsection{Acceptance criterion}

Once a new solution is reconstructed, we either accept it as the new current solution or reject it and reuse the previous one. We use a simulated annealing (SA) based criterion to take this decision. Such an acceptance criterion has been widely used for ALNS studies (see e.g.\cite{ropke2006a} and \cite{iris2017a}).
We accept the new solution $x'$ over the current one $x$ if it is better ($f(x') < f(x)$), or if it is worse with a probability $e^{\frac{-(f(x')-f(x))}{T}}$, where $T$ is the current temperature at a particular iteration, and $f(x)$ is the objective function.
We define an starting an ending temperature,  $T_{start}$ and $T_{end}$ respectively, and the cooling time $t_{cool}$ that defines the duration of going from $T_{start}$ to $T_{end}$.
Based on these parameters, we can define the cooling factor $\tau$ ($0 < \tau < 1$), by isolating it from the formula $T_{end} = T_{start}\tau^{t_{cool}}$. 
This cooling factor allows us to compute the temperature at any given instant. 
Given temperature $T$ at iteration $i$, we find the temperature $T'$ to be used at iteration $i+1$ by computing $T' = T \tau^{t_{it}}$, where $t_{it}$ is the duration of the iteration $i$. 
Following the strategy used in \cite{iris2017a}, we compute $T_{start}$ and $T_{end}$ based on the cost of the initial solution $f(x_0)$ described in Section \ref{sec:iniSol}, where $\xi$ and $\phi$ define the percentage of the cost used to compute $T_{start} = \xi f(x_0)$ and $T_{end} = \phi f(x_0)$. 

\subsection{Adaptive weight adjustment}

One of the main differences between the ALNS method and the standard Large Neighborhood Search (LNS) is the adaptive component of the former. The performance of the employed removal and insertion operators is measured at each iteration. These measures are then used to update the weight and, therefore, the probability of choosing the respective methods.
The most common way of measuring the performance of a method is to give it a different score depending on the quality of the solution. In our case, we define four reward categories as shown in Table \ref{tab:rewardCat}.
\begin{table}[]
    \centering
    \begin{tabular}{c|c}
        \textbf{Category} & \textbf{Parameter} \\ \hline
        Current best solution & $\psi_1$  \\
        Better than current solution & $\psi_2$ \\
        Not better but accepted solution & $\psi_3$ \\
        Rejected solution & $\psi_4$
    \end{tabular}
    \caption{Method reward categories}
    \label{tab:rewardCat}
\end{table}

Let $R$ and $D$ denote the set of insertion and removal operators.
Each removal and insertion method has a probability $\pi^R_i, \pi^D_i$ respectively of being selected at each iteration. Throughout the algorithm run, the probability of selecting these methods gets updated depending on their performance. In our study, we update the probabilities after a $\Delta_{update}$ time interval. During these iterations we accumulate the sum of $\psi^R_i, \psi^D_i$ rewards for each method, and update the weight $\omega^R_i, \omega^D_i$ of each method as indicated in Equation \ref{eq:operatorUpdate}
\begin{equation}\label{eq:operatorUpdate}
    \omega^R_i = (1 - \lambda)\omega^R_i + \lambda \psi^R_i, \quad \omega^D_i = (1 - \lambda)\omega^D_i + \lambda \psi^D_i
\end{equation}
where $\lambda$ is a parameter between 0 and 1 that denotes the degree of adaptability of the method. If $\lambda=0$, the weight remains equal to the previous one. This means that each method would have the same probability throughout the entire algorithm run, behaving like an LNS with multiple neighborhoods. If $\lambda=1$, the new operator's probability solely depends on the score achieved during the last $\Delta_{update}$ and not on previous scores. It is common to use an intermediate value for $\lambda$ strictly between 0 and 1.
Once the weights are updated, the probability of each repair method $\pi_i^R$ and destroy method $\pi_i^D$ can be computed as indicated in Equation (\ref{eq:probUpdate}).
\begin{equation}\label{eq:probUpdate}
    \pi^R_i = \frac{\omega_i}{\sum_{i \in R} \omega_i}, \quad \pi^D_i = \frac{\omega_i}{\sum_{i \in D} \omega_i}
\end{equation} 

\subsection{Local search}\label{sec:localsearch}

An extension of the method is implemented where we perform a local search procedure after reconstructing a new solution. This step aims to incrementally improve the solution by testing small adjustments to the port visits.
\begin{figure}
     \centering
     \begin{subfigure}[b]{0.45\textwidth}
         \centering
         \includegraphics[width=\textwidth]{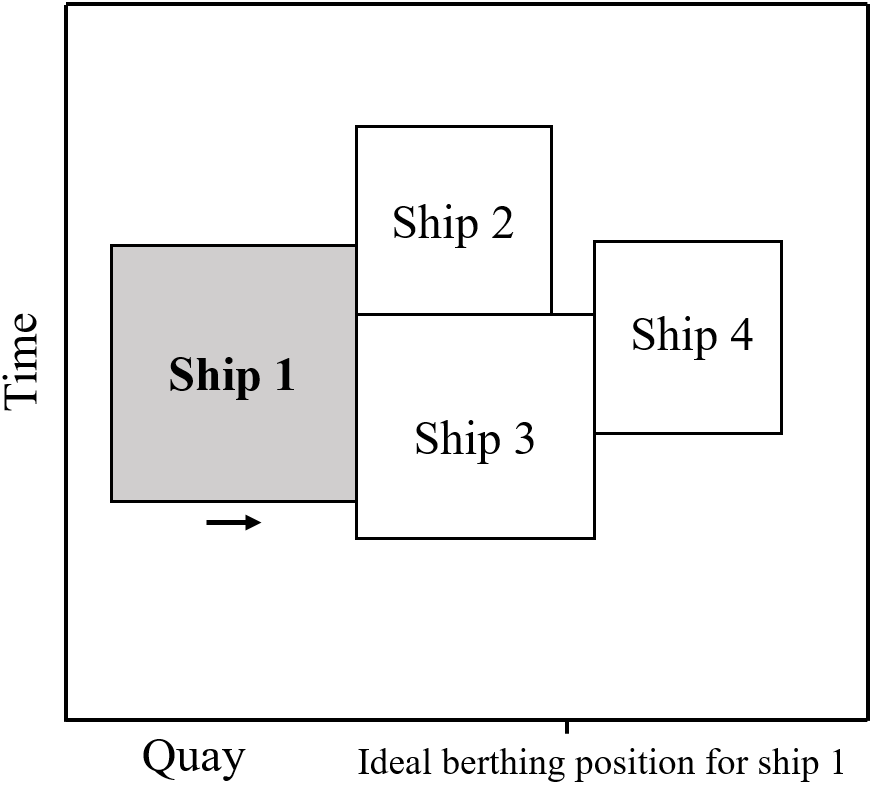}
         \caption{The direction of move for ship 1 is defined.}
         \label{fig:LSex1}
     \end{subfigure}
     \hfill
     \begin{subfigure}[b]{0.45\textwidth}
         \centering
         \includegraphics[width=\textwidth]{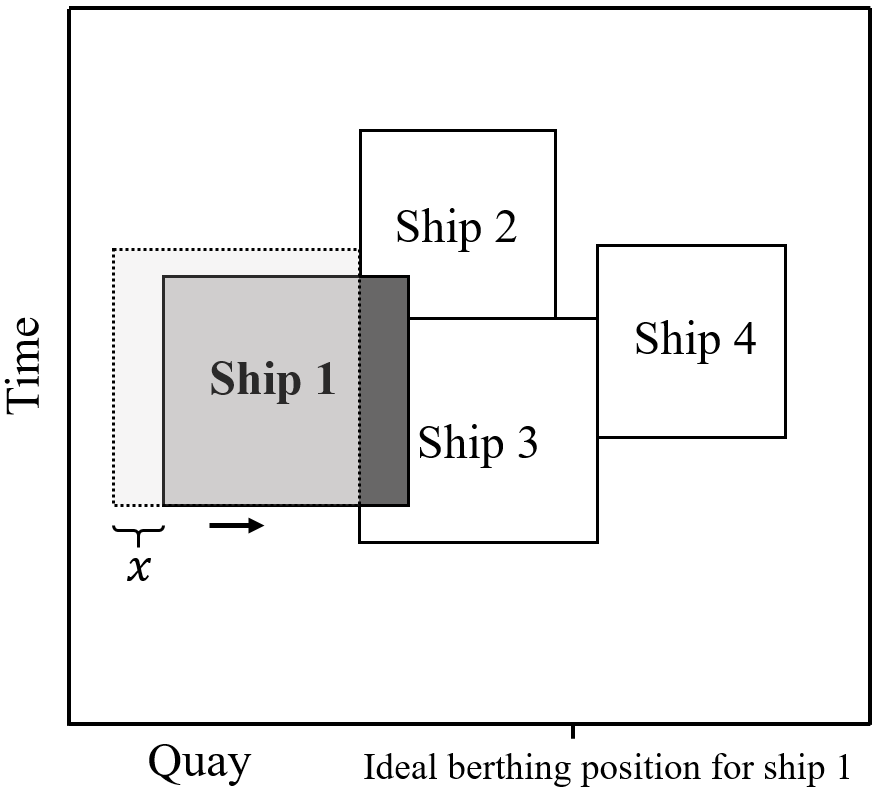}
         \caption{The move generates conflicts with ships 2 and 3.}
         \label{fig:LSex2}
     \end{subfigure} \\
     \begin{subfigure}[b]{0.45\textwidth}
         \centering
         \includegraphics[width=\textwidth]{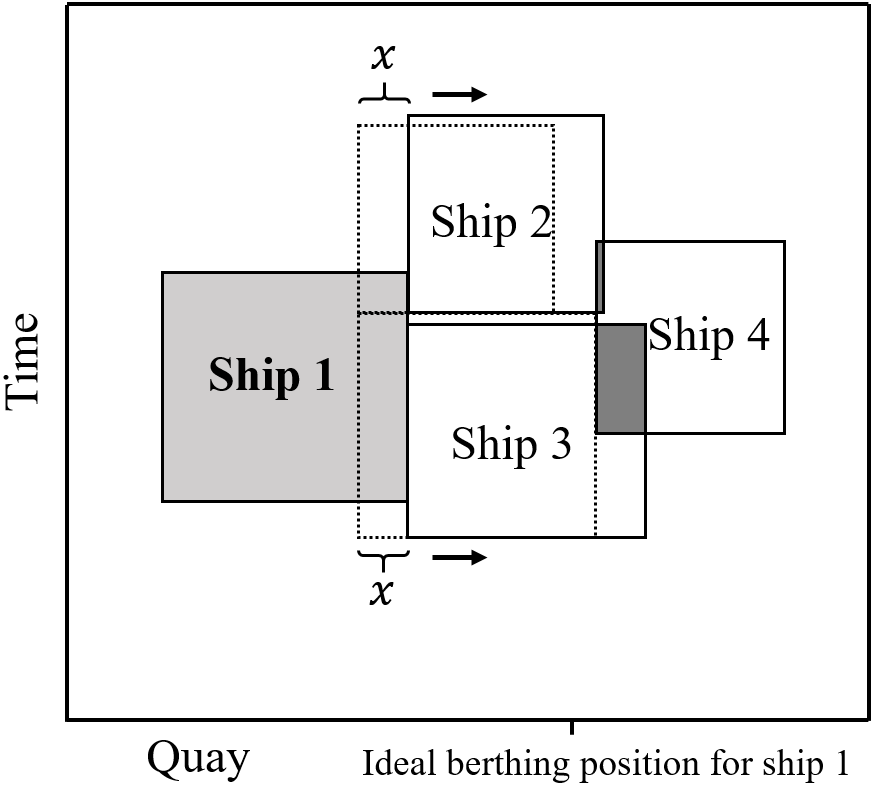}
         \caption{Moving ships 2 and 3 in the same direction generates a new conflict with ship 4.}
         \label{fig:LSex3}
     \end{subfigure}
     \hfill
     \begin{subfigure}[b]{0.45\textwidth}
         \centering
         \includegraphics[width=\textwidth]{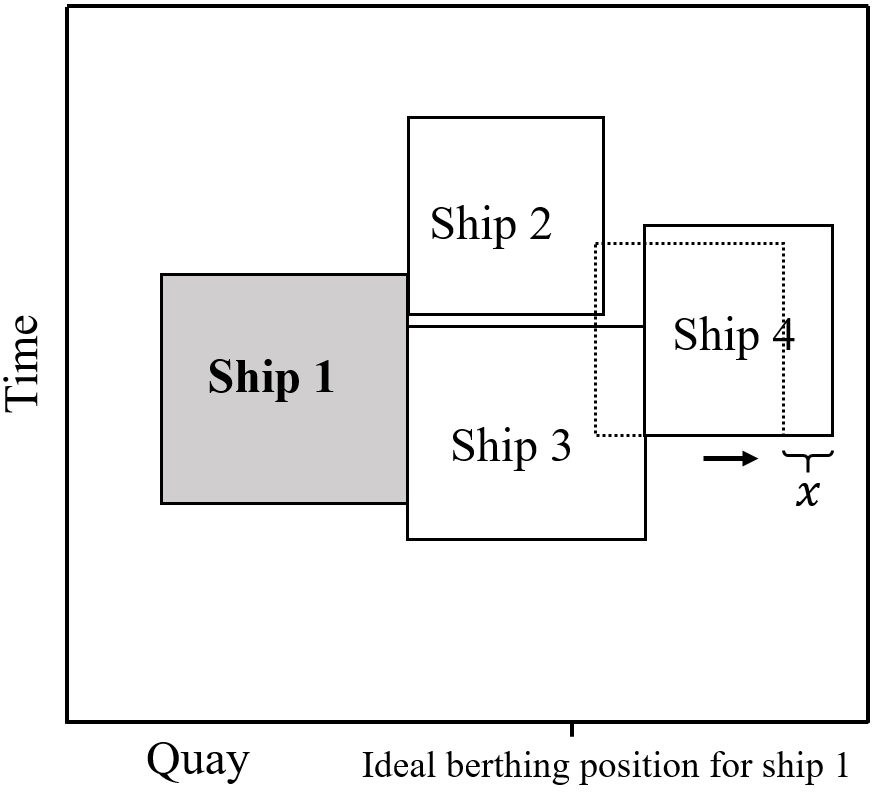}
         \caption{Moving ship 4 in the same direction ends the ejection chain by finding a new feasible solution.}
         \label{fig:LSex4}
     \end{subfigure}
\caption{Example representation of a local search step. The chain of moves originates from ship 1 being moved one step $x$ towards its ideal berthing position. The port visit in grey depicts the first ship to move, and the dark grey indicates an overlapping area. The dashed rectangles represent the original position of the ship before the move.}
\label{fig:LocalSearchEx}
\end{figure}
The procedure is based on the \textit{ejection chains} strategy used in many routing and network-based problems (see \cite{glover1992a, rego1998a, braeysy2003a}). The idea, in our case, is to perturbate the solution by re-planning a port visit to a better position (i.e., lower operational cost) and iteratively re-plan any port visits that conflict with the change. The chain of perturbations is limited to a maximum number of port visits to re-plan $K_{chain}$, and it terminates if this limit is reached or if a conflict-free solution is achieved. Figure \ref{fig:LocalSearchEx} shows an example of this move. Note that the handling time (i.e., the vertical dimension of the port visitsships) is reduced or increased for the ships as their position changes with respect to their ideal position.
A pseudo-code of the procedure is described in Algorithm \ref{alg:localSearch}.
\begin{algorithm}
\DontPrintSemicolon
\KwData{$sol, K_{chain}$: current solution, and the length of the ejection chain (i.e., the maximum number of port visit moves)}
\KwResult{$sol$: resulting solution}
\Begin{
$done \gets false$ \tcp*{initialize termination criterion}
\While(){not $done$}{
    $nextSol \gets sol$ \tcp*{initialize current best solution}
    $\Delta \gets 0$ \tcp*{initialize delta cost variation}
    \For(){$p \in $P}{
        \For(){$n \in N_p$}{
            $toMove = [(p,n)]$ \tcp*{track the port visits to re-plan}
            $sol' \gets sol$ \tcp*{initialize copy of current solution}
            \While{$k \leq K_{chain}$ and $toMove \neq \emptyset$}{
                $(p,n) \gets pop(toMove)$ \tcp*{ get port visit to re-plan}
                $sol' \gets movePortVisit(p,n, sol')$ \tcp*{move port visit}
                $toMove \gets computeConflicts(sol')$ \tcp*{check for conflicts}
                $k \gets k + |toMove|$ \tcp*{update the ejection chain length}
            }
            \If(){$toMove = \emptyset$}{
                $\delta \gets computeDeltaCost(sol,sol')$ \tcp*{compute cost variation}
                \If(){$\delta < \Delta$}{
                    $\Delta \gets \delta$ \tcp*{update best delta cost}
                    $nextSol \gets sol'$ \tcp*{update current best solution}
                }
            }
        }
    }
    \uIf(){$\Delta < 0$}{
        $sol \gets nextSol$ \tcp*{update solution to return}
    }
    \Else(){
        $done \gets true$ \tcp*{no improving neighbor solution}
    }
}
}
\caption{Local search procedure} \label{alg:localSearch}
\end{algorithm}
The function $movePortVisit(p,n)$ performs the perturbation for a given port visit (i.e., ship $n$ at port $p$). It should be noted that the direction is given by the first perturbation made. To find the direction of the first perturbation, we compute the cost variation of moving the port visit in three directions: (i) one segment length towards the ideal position along the spatial axis, and 
(ii) one time instant earlier and (iii) one time instant later along the temporal axis. The direction in the spatial axis is checked if the port visit is not scheduled already at its ideal position. Once the perturbation is performed in the chosen direction, the following port visits in conflict are perturbed in the same direction.

A high value of $K_{chain}$ increases the probability of finding a better solution and the number of operations to compute. The parameter $K_{chain}$ should leverage both solution quality and low computational complexity. Therefore, we define the value of $K_{chain}$ to depend on the number of instances ships and equal to $K_{chain} = 2 \cdot |N|$. The reason for $K_{chain} > |N|$ is that for some movements, a conflicting port visit may require multiple perturbations to achieve a feasible new position, and selecting a lower $K_{chain}$ value may be too restrictive.

\subsection{Algorithm overview}

The overview of the solution method is summarized in Algorithm \ref{alg:alns}. Due to the additional computational effort of the local search procedure, we do not execute it at each iteration. Instead, we only perform it if the reconstructed solution is better than the current one. This reduces the number of times that the local search is performed, allowing the algorithm to perform more iterations while at the same time filtering the times the local search is performed to those where we already have promising solutions.

\begin{algorithm}
\DontPrintSemicolon
\KwData{$inst, param$: a problem instance and a parameter setting for the algorithm}
\KwResult{$bestSol$: best found solution}
\Begin{
$\psi, \pi \gets initialize(inst)$ \tcp*{initialize operator selection parameters}
$sol \gets constructHeuristic(inst)$ \tcp*{construct initial solution}
$bestSol \gets sol$ \\
\While(){timelimit not reached }{
    $currSol \gets sol$ \\
    $removal, insertion \gets selectOperator(\pi)$ \tcp*{select operators}
    $sol \gets insertion(removal(currSol))$ \tcp*{get new solution}
    \If(){$c(sol) < c(currSol)$}{
        $sol \gets localSearch(sol)$
    }
    \If(){isAccepted(sol)}{
        \If(){$c(sol) < c(bestSol)$}{
            $bestSol \gets sol$
        }
        $currSol \gets sol$
    }
    $\pi \gets updateOperatorParams(\psi)$
}
}
\caption{Adaptive large neighborhood search procedure} \label{alg:alns}
\end{algorithm}

\section{Computational results}\label{sec:results}

In this section, we first describe the generation process for the set of benchmark instances, and we then perform a computational study
where we cover both the performance of the method and practical insights of the problem.

\subsection{Instance generation}\label{sec:instGen}
To the best of our knowledge, \cite{martin-iradi22b} is the only study on the MCBAP. The instances presented in the study are rather small and limited. Therefore, we develop a more comprehensive set of benchmark instances. In the absence of real-life data, one could extend current benchmarks instances of the continuous BAP to multiple ports. Instead, we decided to use the public access to port data \citep{marineTraffic}, and we validated it with additional data from an industrial research partner to create an instance generator for the MBAP with a continuous quay.

We consider three different ship types:
(i) feeders or small vessels with a length of up to 200 meters, (ii) medium-size vessels with a length between 200 and 300 meters, and (iii) large vessels longer than 300 meters.
Each ship type has a different speed-fuel consumption relation.
Moreover, we consider three terminals at the three main ports in the north sea:
(i) \textit{Rotterdam APMT} with a quay length of 1600 meters \citep{rotterdam}, (ii) Bremerhaven NTB, with a quay length of 1800 meters \citep{bremen}, and Hamburg EGH, with a quay measuring 2100 meters \citep{hamburg}.
These three ports are relatively close to each other, and large, medium, and small vessels visit them in different sequences as part of their routes.

The duration of the vessel time window is based on the planned port call duration. The planned duration of a vessel's port call can often be updated the days previous to the arrival time. Therefore, we establish a fixed point in time for each ship two weeks before the actual arrival time and retrieve the planned port call duration as the time difference between 
the estimated time of arrival (ETA) and 
the estimated time of departure (ETD). 
We compute this by averaging the planned port call duration for each port and ship type berthing in a period of three months (January-March 2021).
This value is also used to define the minimum handling time $h_0$ (see Table \ref{tab:minH}). 
Port service times that exceed 48, 72, and 96 hours for small, medium, and large ship types, respectively, are categorized as outliers and removed from the dataset. The reason for this is that such long service times usually involve maintenance or fueling operations that are not usually performed on a regular basis. Thus, they are not part of the problem.
\begin{table}[th!]
    \centering
    \caption{Minimum handling type in hours per ship type and terminal. These values define $h^{i,c}_0$.}
    \label{tab:minH}
    \begin{tabular}{c|ccc}
        \textbf{Ship type \textbackslash Terminal} & \textbf{DEHAM} & \textbf{DEBRV} & \textbf{NLRTM} \\ \hline
        \textbf{Feeder} & 10.1 & 12.1 & 10.4 \\
        \textbf{Medium} & 18.0 & 21.8 & 18.4 \\
        \textbf{Large} & 41.0 & 33.7 & 26.7 \\
    \end{tabular}
\end{table}

We define six different ship patterns, each with a given route,
type of ship, and length. All ships visit two or three ports in different orders.
The $N$ ships for a given instance are sampled from the six patterns.

For each ship, we randomize the (i) desired berthing position at each port visited and (ii) the earliest start time $EST_i^c$, following parameters ensuring that feasible sailing times between ports exist.
The estimated finish time $EFT_i^c$ is computed by adding the corresponding value from Table \ref{tab:minH} to $EST_i^c$ ($EFT_i^c = EST_i^c + h_0^{i,c}$).
In addition, we compute the latest finish time $LFT_i^c$ by adding half of the difference between the minimum $h_0^{i,c}$ and maximum $h_{max}^{i,c}$ handling time that the ship can take at the port to $EFT_i^c$ ($LFT_i^c = EFT_i^c + \frac{h_{max}^{i,c} - h_0^{i,c}}{2}$).

As input to the instance generator, we define the number of external ships at each port $Nout$ that are considered fixed. For each external ship, we randomly define: (i) the berthing position and time, (ii) the length (comprised between 180 and 330 meters), and (iii) the handling time, adapted to be proportional to the length of the vessel.

We assume the entire quay is available for berthing unless an external ship is occupying it.
We also consider 10 different speed levels, ranging uniformly between 17-21.5 knots. The ALNS heuristic we present can handle any continuous value of the speed but not the MIP formulation we present. To ensure a fair comparison of the methods we employ a discretized set of speed levels in both the formulation and the solution method.
Furthermore, the distance between ports is computed based on the actual sea distance of the routes.

\subsubsection{Handling time}
It is a general practice, especially on the discrete version of the BAP, to define a different handling time depending on the berthing position. For the continuous version implemented in this paper, we follow the handling time definition presented in \cite{meisel2009a} where deviations from a preferred berthing position are penalized using a deviation factor $\beta \geq 0$ (relative increase in handling time per unit of distance, i.e., meters).
Given the minimum handling time $h^{i,c}_0$ at the preferred berthing position and the actual deviation from the chosen position $\Delta b$ (measured in meters), the handling time is computed as follows:
\begin{equation}
    h^c_i = (1 + \beta \Delta b) h^{i,c}_0
\end{equation}

As a reference, the handling time ranges between 20 and 60 hours for medium vessels and between 30 and 110 hours for large ones. This is given by berthing at the best and worst places, respectively.

\subsubsection{Time windows}

In the MCBAP, there are two types of time windows:
\begin{itemize}
    \item \textit{The time window for each ship at each visited port.} This is given by the port call duration. The time window start must be respected, but the end can be exceeded. The berthing period can, therefore, exceed the end of the time window counting the additional time as a delay.
        \subitem We also consider a time window end that defines the latest finish time (LFT). Ideally, this time window must be respected. However, we allow violating this time window by adding a very high penalty cost. 
    \item \textit{The time window of the berthing positions.} This time window can be seen as the operational hours of a given berthing position or segment of the quay. We assume that all berthing positions are available at any time. It should be noted that potential maintenance windows or partial closures of the quay can be modeled in the same way as an external ship occupying the given positions and time period.
\end{itemize}

\subsubsection{Benchmark instances}
Using the instance generator described in this section, we create a set of benchmark instances.
The entire set comprises 240 instances. Each instance is a combination of the parameter values listed in Table \ref{tab:instParams}: (i) a randomized seed, (ii) the number of ships to optimize, (iii) and the number of external fixed port visits per port, and (iv) the length of the quay segment used to define possible berthing positions.
\begin{table}[]
\centering
\caption{Parameter settings of the benchmark instance set.}
\label{tab:instParams}
\begin{tabular}{c|c|c|c|c}
\textbf{Parameter} & \textbf{Seed} & \textbf{\begin{tabular}[c]{@{}c@{}}Number \\ of ships\end{tabular}} & \textbf{\begin{tabular}[c]{@{}c@{}}Number of external \\ ships per port\end{tabular}} & \textbf{\begin{tabular}[c]{@{}c@{}}Distance between \\ positions\end{tabular}} \\ \hline
\textbf{Values} & 1-10 & 30, 50, 70 & 5, 10 & 10, 20, 40, 80
\end{tabular}%
\end{table}

\subsection{Parameter tuning}

The ALNS algorithm has a total of 18 algorithm parameters. These include parameters used to calibrate the operators, the selection and acceptance criteria, and the weight of the scores for new solutions (see Table \ref{tab:paramtun}).
To select the best value setting for the algorithm parameters we conducted a parameter tuning.
We selected a subset of 12 instances that are representative of the entire set. For the tuning, we run the automatic algorithm configurator \textit{Pydgga} \citep{Ansotegui2021} for the 30 generations with a time limit of 5 minutes for each ALNS run. The configurator allows to parallelize the process, and the entire parameter tuning lasted 8 hours. In Table \ref{tab:paramtun}, we define the domain of each algorithm parameter used for tuning and the found setting.
The models and solution methods are written in \textit{Julia} and run in a 2.90
GHz Intel Xeon Gold 6226R using one thread and 16 GB of RAM.
\begin{table}[]
\centering
\caption{Studied range and chosen setting of algorithm parameters after the parameter tuning.}
\label{tab:paramtun}
\resizebox{\textwidth}{!}{%
\begin{tabular}{ll|rr|r}
\textbf{Symbol} & \textbf{Description} & \textbf{Min val} & \textbf{Max val} & \textbf{Tuned setting} \\ \hline
 $\epsilon$ & value used to compute the cooling ratio & 0.005 & 0.2 & \textbf{0.157} \\
 $\varphi$ & \begin{tabular}[c]{@{}l@{}}pct of initial solution obj used to \\ define start temperature (when reheated)\end{tabular} & 0.01 & 0.05 & \textbf{0.0246} \\
 $\xi$ & \begin{tabular}[c]{@{}l@{}}pct of initial solution obj used to\\ define end temperature (to be reheated)\end{tabular} & 0.00005 & 0.001 & \textbf{0.000269} \\
 $\rho$ & \begin{tabular}[c]{@{}l@{}}degree of destruction, pct of total port\\ visits to be removed by the removal methods\end{tabular} & 0.3 & 0.6 & \textbf{0.326} \\
 $A$& weight for position deviation in shaw removal & 0.5 & 2 & \textbf{0.55} \\
 $B$ & \begin{tabular}[c]{@{}l@{}}weight for berthing start time deviation\\ in shaw removal\end{tabular} & 0.5 & 2 & \textbf{1.36} \\
 $C$ & \begin{tabular}[c]{@{}l@{}}weight for berthing end time deviation in\\ shaw removal\end{tabular} & 0.5 & 2 & \textbf{0.89} \\
 $\alpha$ & \begin{tabular}[c]{@{}l@{}}randomness parameter for shaw removal \\ method\end{tabular} & 1 & 3 & \textbf{2.66} \\
 $\gamma$ & \begin{tabular}[c]{@{}l@{}}randomness parameter for random greedy\\ repair method\end{tabular} & 1 & 3 & \textbf{2.85} \\
 $\mu$ & \begin{tabular}[c]{@{}l@{}}randomness parameter for arrival greedy\\ repair method\end{tabular} & 1 & 3 & \textbf{2.6} \\
 $\kappa$ & k-regret parameter & 2 & 4 & \textbf{2} \\
 $\Delta$ & \begin{tabular}[c]{@{}l@{}}number of iterations between updating the \\ weights of each method\end{tabular} & 0.01 & 0.05 & \textbf{0.017} \\
 $\eta$ & \begin{tabular}[c]{@{}l@{}}parameter to adjust the importance of \\ recent scores vs. previous weight\end{tabular} & 0.3 & 0.7 & \textbf{0.456} \\
 $\psi_1$ & score when finding a current best solution & 10 & 20 & \textbf{11} \\
 $\psi_2$ & \begin{tabular}[c]{@{}l@{}}score when finding a solution better than\\ the current solution\end{tabular} & 4 & 8 & \textbf{4} \\
 $\psi_3$ & score when the solution is accepted & 1 & 3 & \textbf{2} \\
 $\psi_4$ & score when the solution is rejected & - & - & \textbf{0} \\
 $\beta$ & \begin{tabular}[c]{@{}l@{}}parameter that defines the position\\ bounds for ship (times the length)\end{tabular} & 2 & 5 & \textbf{4.02}
\end{tabular}}
\end{table}

\subsection{Method performance}

To measure the quality of the method, we compare the presented method with its variants and with a baseline commercial solvers; \textit{CPLEX v12.10}.

Table \ref{tab:alnsVSmip} compares the results between CPLEX and the ALNS method.
We compute the objective gap for each instance run as $\frac{z_{obj} - z_{best}}{z_{best}}$ where $z_{obj}$ is the objective value of the best solution of the run, and $z_{best}$ is the best-known solution across all experiments.
We have grouped all instances per number of ships, external ships, and distance between berthing positions. Each group contains 10 instances and is named \textit{X-Y-Z} according to their common characteristics. \textit{X} is the number of ships, \textit{Y} is the number of external ships per port, and \textit{Z} is the distance between consecutive positions considered in meters (i.e., segment length).
Therefore, each row in the table corresponds to the average value across instances with different seed values. The ALNS is tested by running each instance 10 times and computing the average, best and worst run.

\begin{table}[]
\centering
\caption{Gap for the MIP formulation solved by CPLEX and the ALNS method with a time limit of 5 minutes and 1 hour. We also report the average gap between the best and worst runs of the ALNS. Each row corresponds to an instance group (i.e., the average results across 10 instances of the same size). 
}
\label{tab:alnsVSmip}
\resizebox{\textwidth}{!}{%
\begin{tabular}{c|rrrr|rrrr}
\textbf{} & \multicolumn{4}{c|}{\textbf{5 minutes}} & \multicolumn{4}{c}{\textbf{1 hour}} \\ \hline
\textbf{\begin{tabular}[c]{@{}c@{}}Instance\\ group\end{tabular}} & \multicolumn{1}{c}{\textbf{\begin{tabular}[c]{@{}c@{}}MIP\\ gap (\%)\end{tabular}}} & \multicolumn{1}{c}{\textbf{\begin{tabular}[c]{@{}c@{}}ALNS\\ gap (\%)\end{tabular}}} & \multicolumn{1}{c}{\textbf{\begin{tabular}[c]{@{}c@{}}Best ALNS\\ gap (\%)\end{tabular}}} & \multicolumn{1}{c|}{\textbf{\begin{tabular}[c]{@{}c@{}}Worst ALNS\\ gap (\%)\end{tabular}}} & \multicolumn{1}{c}{\textbf{\begin{tabular}[c]{@{}c@{}}MIP\\ gap (\%)\end{tabular}}} & \multicolumn{1}{c}{\textbf{\begin{tabular}[c]{@{}c@{}}ALNS\\ gap (\%)\end{tabular}}} & \multicolumn{1}{c}{\textbf{\begin{tabular}[c]{@{}c@{}}Best ALNS\\ gap (\%)\end{tabular}}} & \multicolumn{1}{c}{\textbf{\begin{tabular}[c]{@{}c@{}}Worst ALNS\\ gap (\%)\end{tabular}}} \\ \hline
30\_5\_10 & \textbf{10.1} & 11.9 & 7.3 & 17.4 & \textbf{1.5} & 4.8 & 1.8 & 8.4 \\
30\_5\_20 & \textbf{7.4} & 9.4 & 4.6 & 14.9 & \textbf{1.7} & 3.3 & 0.9 & 6.0 \\
30\_5\_40 & 10.2 & \textbf{6.9} & 3.6 & 10.2 & \textbf{1.2} & 3.4 & 1.2 & 5.5 \\
30\_5\_80 & 15.3 & \textbf{7.5} & 4.3 & 10.5 & \textbf{2.1} & 3.4 & 0.9 & 5.5 \\
30\_10\_10 & 15.1 & \textbf{10.0} & 6.1 & 14.5 & \textbf{2.7} & 3.6 & 0.8 & 6.8 \\
30\_10\_20 & 16.5 & \textbf{8.1} & 4.7 & 11.9 & 5.8 & \textbf{2.9} & 0.5 & 5.2 \\
30\_10\_40 & 13.3 & \textbf{5.9} & 3.3 & 8.4 & 3.5 & \textbf{2.3} & 0.2 & 4.3 \\
30\_10\_80 & 16.9 & \textbf{5.2} & 3.3 & 7.6 & 5.7 & \textbf{2.3} & 0.5 & 4.5 \\
50\_5\_10 & 22.5 & \textbf{14.4} & 8.6 & 19.3 & \textbf{2.2} & 4.3 & 0.7 & 9.6 \\
50\_5\_20 & 45.7 & \textbf{11.4} & 6.1 & 18.2 & 6.4 & \textbf{3.6} & 0.8 & 6.9 \\
50\_5\_40 & 59.1 & \textbf{8.8} & 4.3 & 13.1 & 11.1 & \textbf{2.9} & 0.0 & 6.2 \\
50\_5\_80 & 82.7 & \textbf{7.2} & 4.3 & 10.6 & 14.3 & \textbf{2.9} & 0.5 & 5.3 \\
50\_10\_10 & 66.5 & \textbf{10.3} & 5.2 & 14.6 & 8.1 & \textbf{2.9} & 0.6 & 5.9 \\
50\_10\_20 & 61.4 & \textbf{8.3} & 4.4 & 13.0 & 9.2 & \textbf{2.2} & 0.0 & 4.6 \\
50\_10\_40 & 74.6 & \textbf{7.0} & 3.9 & 11.3 & 13.4 & \textbf{2.7} & 0.0 & 5.5 \\
50\_10\_80 & 92.3 & \textbf{6.2} & 3.5 & 9.3 & 16.7 & \textbf{2.2} & 0.0 & 4.3 \\
70\_5\_10 & 139.5 & \textbf{13.8} & 7.2 & 18.9 & 8.7 & \textbf{3.9} & 0.8 & 8.2 \\
70\_5\_20 & 103.4 & \textbf{10.9} & 5.6 & 16.3 & 9.3 & \textbf{2.4} & 0.0 & 5.9 \\
70\_5\_40 & 141.3 & \textbf{7.8} & 4.2 & 13.4 & 14.7 & \textbf{2.4} & 0.0 & 4.8 \\
70\_5\_80 & 136.5 & \textbf{5.8} & 2.9 & 9.0 & 17.9 & \textbf{2.3} & 0.0 & 4.5 \\
70\_10\_10 & 95.2 & \textbf{11.7} & 7.6 & 16.0 & 15.5 & \textbf{2.7} & 0.0 & 5.5 \\
70\_10\_20 & 90.8 & \textbf{8.8} & 4.2 & 13.0 & 16.8 & \textbf{2.3} & 0.0 & 5.1 \\
70\_10\_40 & 122.0 & \textbf{7.3} & 3.9 & 11.1 & 29.5 & \textbf{2.3} & 0.0 & 4.7 \\
70\_10\_80 & 141.0 & \textbf{5.5} & 3.3 & 8.4 & 28.2 & \textbf{2.1} & 0.0 & 3.9 \\ \hline
\textbf{Average} & 65.8 & \textbf{8.7} & 4.8 & 13.0 & 10.3 & \textbf{2.9} & 0.4 & 5.7
\end{tabular}%
}
\end{table}

We observe that CPLEX scales poorly, especially in instances with more than 30 ships. The gap is better for the smallest instances but quickly worsens. On average, the ALNS method outperforms the commercial solver by achieving tighter gaps in most instances. 
The gap is an indicator of the method performance relative to each other but does not provide an optimality guarantee. The lower bounds obtained with CPLEX indicate a high optimality gap. This could be due to a low-quality solution or a poor lower bound.
\cite{martin-iradi2022a} indicated that the relaxation of the MIP formulation for the MBAP with a discrete quay could be poor and showed that a network-flow reformulation could tighten the relaxation significantly.
As indicated in \cite{martin-iradi22b}, network-flow formulations for the MCBAP can suffer from scalability but show that the relaxation is stronger.
We compare the results from the branch-and-price method presented in \cite{martin-iradi22b} with the ALNS method. 
The formulation presented in \cite{martin-iradi22b} is slightly different from the one addressed in this paper. 
The formulation from \cite{martin-iradi22b} defines the latest finish time for each ship berthing at a port that must be satisfied. 
We have adapted the method from \cite{martin-iradi22b} to the formulation of this study.
The branch-and-price method is based on a graph representation, and therefore, we need to establish the latest possible berthing time. This is set to 50 \% more than the latest finish time. This allows the method to exceed the latest finish time while maintaining the graph at a reasonable size. This is a generous bound, and our empirical studies show that this bound does not affect the optimal solution. 
For that, we have also generated a new set of instances of similar size to the ones presented in \cite{martin-iradi22b} using the instance generator defined in Section \ref{sec:instGen}. 
\begin{table}[]
\centering
\caption{Parameter settings of the instance set based on the ones from \cite{martin-iradi22b}.}
\label{tab:instBP}
\begin{tabular}{c|c|c|c|c}
\textbf{Parameter} & \textbf{Seed} & \textbf{\begin{tabular}[c]{@{}c@{}}Number \\ of ships\end{tabular}} & \textbf{\begin{tabular}[c]{@{}c@{}}Number of external \\ ships per port\end{tabular}} & \textbf{\begin{tabular}[c]{@{}c@{}}Distance between \\ positions\end{tabular}} \\ \hline
\textbf{Values} & 1-5 & 4-15 & 3-5 & 10, 20, 40, 80
\end{tabular}%
\end{table}
The input parameters for the instance set are defined in Table \ref{tab:instBP}. The entire set comprises 720 instances, one for each combination of input parameters.
The results are shown in Table \ref{tab:alnsVsBP}, where we have grouped the instances in batches of 60 according to the number of ships. We compare the branch-and-price method with the ALNS and MIP formulation presented in this study.
For both the branch-and-price and CPLEX we compute their optimality gap, where we observe that the branch-and-price method achieves a better gap due to the tighter lower bound in most cases. 
We also compute the gap to the best-known solution for all three methods.
In this case, we observe that CPLEX provides the best performance, showing that despite its poorer relaxation, the upper bounds found are near-optimal.
The branch-and-price method still shows a robust performance but for short computational times and larger instances, the ALNS method is able to provide better solutions.

\begin{table}[]
\centering
\caption{Performance comparison across 720 instances between the MIP formulation solved by CPLEX, the adapted branch-and-price method from \cite{martin-iradi22b}, and the ALNS method, with a time limit of 5 minutes and 1 hour. Each row shows the average gap values across all instances with same number of ships.}
\label{tab:alnsVsBP}
\resizebox{\textwidth}{!}{%
\begin{tabular}{r|rr|rr|rr|rr|rr}
\multicolumn{1}{c|}{\multirow{2}{*}{\textbf{\begin{tabular}[c]{@{}c@{}}Number \\ of ships\end{tabular}}}} & \multicolumn{2}{c|}{\textbf{\begin{tabular}[c]{@{}c@{}}CPLEX\\ Opt. gap (\%)\end{tabular}}} & \multicolumn{2}{c|}{\textbf{\begin{tabular}[c]{@{}c@{}}Branch-and-price\\ Opt. gap (\%)\end{tabular}}} & \multicolumn{2}{c|}{\textbf{\begin{tabular}[c]{@{}c@{}}CPLEX\\ Gap (\%)\end{tabular}}} & \multicolumn{2}{c|}{\textbf{\begin{tabular}[c]{@{}c@{}}Branch-and-price\\ Gap (\%)\end{tabular}}} & \multicolumn{2}{c}{\textbf{\begin{tabular}[c]{@{}c@{}}ALNS+LS\\ Gap (\%)\end{tabular}}} \\ \cline{2-11} 
\multicolumn{1}{c|}{} & \multicolumn{1}{c|}{5 min.} & \multicolumn{1}{c|}{1 hour} & \multicolumn{1}{c|}{5 min.} & \multicolumn{1}{c|}{1 hour} & \multicolumn{1}{c|}{5 min.} & \multicolumn{1}{c|}{1 hour} & \multicolumn{1}{c|}{5 min.} & \multicolumn{1}{c|}{1 hour} & \multicolumn{1}{c|}{5 min.} & \multicolumn{1}{c}{1 hour} \\ \hline
4 & \multicolumn{1}{r|}{\textbf{0.0}} & \textbf{0.0} & \multicolumn{1}{r|}{\textbf{0.0}} & \textbf{0.0} & \multicolumn{1}{r|}{\textbf{0.0}} & \textbf{0.0} & \multicolumn{1}{r|}{\textbf{0.0}} & \textbf{0.0} & \multicolumn{1}{r|}{0.1} & \textbf{0.0} \\
5 & \multicolumn{1}{r|}{\textbf{0.0}} & \textbf{0.0} & \multicolumn{1}{r|}{\textbf{0.0}} & \textbf{0.0} & \multicolumn{1}{r|}{\textbf{0.0}} & \textbf{0.0} & \multicolumn{1}{r|}{\textbf{0.0}} & \textbf{0.0} & \multicolumn{1}{r|}{\textbf{0.0}} & \textbf{0.0} \\
6 & \multicolumn{1}{r|}{\textbf{0.0}} & \textbf{0.0} & \multicolumn{1}{r|}{\textbf{0.0}} & \textbf{0.0} & \multicolumn{1}{r|}{\textbf{0.0}} & \textbf{0.0} & \multicolumn{1}{r|}{\textbf{0.0}} & \textbf{0.0} & \multicolumn{1}{r|}{0.2} & 0.1 \\
7 & \multicolumn{1}{r|}{\textbf{0.0}} & \textbf{0.0} & \multicolumn{1}{r|}{0.2} & \textbf{0.0} & \multicolumn{1}{r|}{\textbf{0.0}} & \textbf{0.0} & \multicolumn{1}{r|}{\textbf{0.0}} & \textbf{0.0} & \multicolumn{1}{r|}{0.3} & 0.2 \\
8 & \multicolumn{1}{r|}{\textbf{0.0}} & \textbf{0.0} & \multicolumn{1}{r|}{\textbf{0.0}} & \textbf{0.0} & \multicolumn{1}{r|}{\textbf{0.0}} & \textbf{0.0} & \multicolumn{1}{r|}{\textbf{0.0}} & \textbf{0.0} & \multicolumn{1}{r|}{0.3} & 0.2 \\
9 & \multicolumn{1}{r|}{\textbf{0.0}} & \textbf{0.0} & \multicolumn{1}{r|}{0.3} & 0.1 & \multicolumn{1}{r|}{\textbf{0.0}} & \textbf{0.0} & \multicolumn{1}{r|}{0.2} & \textbf{0.0} & \multicolumn{1}{r|}{0.3} & 0.2 \\
10 & \multicolumn{1}{r|}{\textbf{0.0}} & \textbf{0.0} & \multicolumn{1}{r|}{0.1} & \textbf{0.0} & \multicolumn{1}{r|}{\textbf{0.0}} & \textbf{0.0} & \multicolumn{1}{r|}{\textbf{0.0}} & \textbf{0.0} & \multicolumn{1}{r|}{0.5} & 0.3 \\
11 & \multicolumn{1}{r|}{1.0} & 0.2 & \multicolumn{1}{r|}{\textbf{0.3}} & \textbf{0.1} & \multicolumn{1}{r|}{\textbf{0.1}} & \textbf{0.0} & \multicolumn{1}{r|}{\textbf{0.1}} & \textbf{0.0} & \multicolumn{1}{r|}{0.3} & 0.2 \\
12 & \multicolumn{1}{r|}{1.1} & 0.6 & \multicolumn{1}{r|}{\textbf{0.8}} & \textbf{0.2} & \multicolumn{1}{r|}{\textbf{0.1}} & \textbf{0.0} & \multicolumn{1}{r|}{0.4} & \textbf{0.0} & \multicolumn{1}{r|}{0.2} & 0.1 \\
13 & \multicolumn{1}{r|}{4.0} & 2.7 & \multicolumn{1}{r|}{\textbf{1.1}} & \textbf{0.6} & \multicolumn{1}{r|}{\textbf{0.1}} & \textbf{0.1} & \multicolumn{1}{r|}{0.4} & 0.2 & \multicolumn{1}{r|}{0.4} & 0.2 \\
14 & \multicolumn{1}{r|}{3.5} & 2.0 & \multicolumn{1}{r|}{\textbf{1.2}} & \textbf{0.4} & \multicolumn{1}{r|}{\textbf{0.1}} & \textbf{0.1} & \multicolumn{1}{r|}{0.6} & \textbf{0.1} & \multicolumn{1}{r|}{0.7} & 0.4 \\
15 & \multicolumn{1}{r|}{7.3} & 4.8 & \multicolumn{1}{r|}{\textbf{1.7}} & \textbf{1.0} & \multicolumn{1}{r|}{\textbf{0.3}} & \textbf{0.1} & \multicolumn{1}{r|}{0.7} & 0.2 & \multicolumn{1}{r|}{0.5} & 0.3 \\ \hline
Average & \multicolumn{1}{r|}{1.41} & 0.86 & \multicolumn{1}{r|}{\textbf{0.49}} & \textbf{0.20} & \multicolumn{1}{r|}{\textbf{0.06}} & \textbf{0.02} & \multicolumn{1}{r|}{0.20} & 0.05 & \multicolumn{1}{r|}{0.30} & 0.19
\end{tabular}%
}
\end{table}

\begin{table}[]
\centering
\caption{Performance comparison between variants of the proposed ALNS method.}
\label{tab:alnsVars}
\resizebox{\textwidth}{!}{%
\begin{tabular}{c|rrrr|rrrr}
\textbf{} & \multicolumn{4}{c|}{\textbf{5 minutes}} & \multicolumn{4}{c}{\textbf{1 hour}} \\ \hline
\textbf{\begin{tabular}[c]{@{}c@{}}Instance\\ group\end{tabular}} & \multicolumn{1}{c}{\textbf{\begin{tabular}[c]{@{}c@{}}ALNS + LS\\ gap (\%)\end{tabular}}} & \multicolumn{1}{c}{\textbf{\begin{tabular}[c]{@{}c@{}}ALNS (no LS)\\ gap (\%)\end{tabular}}} & \multicolumn{1}{c}{\textbf{\begin{tabular}[c]{@{}c@{}}LNS + LS\\ gap (\%)\end{tabular}}} & \multicolumn{1}{c|}{\textbf{\begin{tabular}[c]{@{}c@{}}LNS (no LS)\\ gap (\%)\end{tabular}}} & \multicolumn{1}{c}{\textbf{\begin{tabular}[c]{@{}c@{}}ALNS + LS\\ gap (\%)\end{tabular}}} & \multicolumn{1}{c}{\textbf{\begin{tabular}[c]{@{}c@{}}ALNS (no LS)\\ gap (\%)\end{tabular}}} & \multicolumn{1}{c}{\textbf{\begin{tabular}[c]{@{}c@{}}LNS + LS\\ gap (\%)\end{tabular}}} & \multicolumn{1}{c}{\textbf{\begin{tabular}[c]{@{}c@{}}LNS (no LS)\\ gap (\%)\end{tabular}}} \\ \hline
30\_5\_10 & 11.9 & \textbf{11.8} & 12.0 & 13.6 & \textbf{4.8} & 4.9 & 5.5 & 5.9 \\
30\_5\_20 & 9.4 & 9.6 & \textbf{8.8} & 10.6 & \textbf{3.3} & 3.6 & 3.4 & 4.0 \\
30\_5\_40 & \textbf{6.9} & 7.3 & 7.3 & 7.8 & 3.4 & 3.3 & \textbf{2.1} & 2.5 \\
30\_5\_80 & 7.5 & 7.5 & \textbf{6.4} & 7.0 & 3.4 & 3.7 & \textbf{2.1} & 2.3 \\
30\_10\_10 & \textbf{10.0} & 10.9 & 10.6 & 11.8 & \textbf{3.6} & 3.7 & 4.2 & 4.7 \\
30\_10\_20 & \textbf{8.1} & 8.2 & 8.6 & 9.6 & 2.9 & \textbf{2.8} & 3.1 & 3.3 \\
30\_10\_40 & 5.9 & 6.3 & \textbf{5.6} & 6.3 & 2.3 & 2.4 & \textbf{1.4} & 1.7 \\
30\_10\_80 & 5.2 & 5.4 & \textbf{4.5} & 5.3 & 2.3 & 2.3 & \textbf{0.8} & 1.0 \\
50\_5\_10 & \textbf{14.4} & 15.5 & 15.1 & 18.2 & \textbf{4.3} & 5.1 & 7.1 & 9.4 \\
50\_5\_20 & 11.4 & \textbf{12.2} & 12.8 & 16.3 & \textbf{3.6} & 4.4 & 5.1 & 7.7 \\
50\_5\_40 & \textbf{8.8} & 10.2 & 11.0 & 13.2 & \textbf{2.9} & 3.9 & 3.8 & 6.2 \\
50\_5\_80 & \textbf{7.2} & 8.3 & 8.2 & 10.6 & \textbf{2.9} & 3.0 & 3.4 & 4.6 \\
50\_10\_10 & \textbf{10.3} & 11.4 & 12.0 & 14.4 & \textbf{2.9} & 3.9 & 5.1 & 7.1 \\
50\_10\_20 & \textbf{8.3} & 9.2 & 9.9 & 12.5 & \textbf{2.2} & 3.3 & 3.6 & 5.7 \\
50\_10\_40 & \textbf{7.0} & 8.1 & 8.1 & 10.9 & \textbf{2.7} & 3.5 & 3.1 & 4.9 \\
50\_10\_80 & \textbf{6.2} & 6.6 & 6.6 & 8.7 & \textbf{2.2} & 2.6 & 2.6 & 3.1 \\
70\_5\_10 & \textbf{13.8} & 14.7 & 14.6 & 17.4 & \textbf{3.9} & 5.4 & 6.7 & 10.9 \\
70\_5\_20 & \textbf{10.9} & 11.4 & 11.4 & 14.4 & \textbf{2.4} & 4.6 & 5.1 & 8.5 \\
70\_5\_40 & \textbf{7.8} & 10.1 & 9.4 & 13.1 & \textbf{2.4} & 3.9 & 3.9 & 7.1 \\
70\_5\_80 & \textbf{5.8} & 7.8 & 7.6 & 11.0 & \textbf{2.3} & 3.5 & 3.2 & 5.5 \\
70\_10\_10 & \textbf{11.7} & 12.6 & 11.8 & 14.8 & \textbf{2.7} & 4.3 & 5.2 & 8.3 \\
70\_10\_20 & \textbf{8.8} & 10.3 & 9.6 & 13.3 & \textbf{2.3} & 4.2 & 4.6 & 7.6 \\
70\_10\_40 & \textbf{7.3} & 9.2 & 8.4 & 11.9 & \textbf{2.3} & 3.8 & 3.9 & 6.2 \\
70\_10\_80 & \textbf{5.5} & 6.8 & 6.6 & 9.4 & \textbf{2.1} & 3.0 & 3.0 & 4.5 \\ \hline
\textbf{Average} & \textbf{8.7} & 9.7 & 9.4 & 11.8 & \textbf{2.9} & 3.7 & 3.8 & 5.5
\end{tabular}%
}
\end{table}
The ALNS method has two main components that differentiate it from other heuristics: (i) the adaptive procedure that guides the operator selection and (ii) the local search procedure that is performed when promising solutions are found. To quantify the impact of these two procedures, we compare the proposed method to its variants with and without each of the components.
One variant is the method without its adaptive component (i.e., large neighborhood search (LNS)), meaning that each removal and insertion operator has an equal probability of being selected throughout the algorithm run. Another variant is the ALNS method without the local search (LS). 
The objective gap across the methods is compared in Table \ref{tab:alnsVars}, with a time limit of 5 minutes, and 1 hour.
For the short time limit, we see that the ALNS without the local search performs the best in most instances. Once the  time limit is increased, the value of the local search is more apparent, and it provides the best results in most instances.

\begin{table}[]
\centering
\caption{Algorithm comparison with different frequencies for the local search procedure with a time limit of one hour.}
\label{tab:LSvars}
\resizebox{\textwidth}{!}{%
\begin{tabular}{c|rr|rr|rr|rrr}
\multicolumn{1}{c|}{\textbf{}} & \multicolumn{2}{c|}{\textbf{\begin{tabular}[c]{@{}c@{}}Local search\\ every iteration\end{tabular}}} & \multicolumn{2}{c|}{\textbf{\begin{tabular}[c]{@{}c@{}}Local search\\ every 2 iterations\end{tabular}}} & \multicolumn{2}{c|}{\textbf{\begin{tabular}[c]{@{}c@{}}Local search\\ every 4 iterations\end{tabular}}} & \multicolumn{3}{c}{\textbf{\begin{tabular}[c]{@{}c@{}}Local search every iteration where\\ the solution is better than the incumbent\end{tabular}}} \\ \hline
\multicolumn{1}{c|}{\textbf{\begin{tabular}[c]{@{}c@{}}Instance\\ group\end{tabular}}} & \multicolumn{1}{c}{\textbf{\begin{tabular}[c]{@{}c@{}}Gap\\ (\%)\end{tabular}}} & \multicolumn{1}{c|}{\textbf{\begin{tabular}[c]{@{}c@{}}Iterations\\ x1000\end{tabular}}} & \multicolumn{1}{c}{\textbf{\begin{tabular}[c]{@{}c@{}}Gap\\ (\%)\end{tabular}}} & \multicolumn{1}{c|}{\textbf{\begin{tabular}[c]{@{}c@{}}Iterations\\ x1000\end{tabular}}} & \multicolumn{1}{c}{\textbf{\begin{tabular}[c]{@{}c@{}}Gap\\ (\%)\end{tabular}}} & \multicolumn{1}{c|}{\textbf{\begin{tabular}[c]{@{}c@{}}Iterations\\ x1000\end{tabular}}} & \multicolumn{1}{c}{\textbf{\begin{tabular}[c]{@{}c@{}}Gap\\ (\%)\end{tabular}}} & \multicolumn{1}{c}{\textbf{\begin{tabular}[c]{@{}c@{}}Iterations\\ x1000\end{tabular}}} & \multicolumn{1}{c}{\textbf{\begin{tabular}[c]{@{}c@{}}\% of  iterations\\ with local search\end{tabular}}} \\ \hline
30\_5\_10 & 5.3 & 5.2 & 5.3 & 7.3 & 5.3 & 7.8 & 4.8 & 12.7 & 1.3 \\
30\_5\_20 & 4.2 & 8.0 & 4.1 & 11.1 & 4.2 & 12.8 & 3.3 & 18.6 & 1.5 \\
30\_5\_40 & 3.8 & 11.0 & 3.4 & 16.5 & 3.6 & 21.1 & 3.4 & 27.4 & 2.5 \\
30\_5\_80 & 4.3 & 15.7 & 3.9 & 25.7 & 3.8 & 36.0 & 3.4 & 60.4 & 2.7 \\
30\_10\_10 & 3.9 & 5.6 & 4.0 & 7.1 & 4.0 & 7.8 & 3.6 & 10.9 & 1.6 \\
30\_10\_20 & 3.1 & 8.8 & 3.1 & 12.3 & 3.5 & 13.3 & 2.9 & 20.0 & 1.7 \\
30\_10\_40 & 2.5 & 12.4 & 2.3 & 18.2 & 2.2 & 23.9 & 2.3 & 30.4 & 2.9 \\
30\_10\_80 & 2.7 & 17.6 & 2.5 & 28.9 & 2.5 & 39.4 & 2.3 & 63.4 & 3.0 \\
50\_5\_10 & 3.9 & 2.0 & 4.4 & 3.2 & 4.9 & 3.8 & 4.3 & 8.5 & 0.4 \\
50\_5\_20 & 4.2 & 2.6 & 4.5 & 4.3 & 4.6 & 6.2 & 3.6 & 16.2 & 0.4 \\
50\_5\_40 & 3.4 & 3.4 & 2.9 & 5.3 & 3.4 & 7.8 & 2.9 & 13.8 & 0.8 \\
50\_5\_80 & 4.1 & 4.5 & 3.9 & 7.6 & 3.7 & 11.6 & 2.9 & 24.8 & 1.1 \\
50\_10\_10 & 2.9 & 2.3 & 2.9 & 3.5 & 3.7 & 4.6 & 2.9 & 8.8 & 0.7 \\
50\_10\_20 & 2.6 & 3.1 & 2.8 & 5.1 & 2.9 & 6.3 & 2.2 & 15.4 & 0.6 \\
50\_10\_40 & 3.1 & 4.2 & 3.1 & 6.4 & 3.1 & 9.1 & 2.7 & 14.3 & 1.2 \\
50\_10\_80 & 2.9 & 5.5 & 2.8 & 9.2 & 3.0 & 13.4 & 2.2 & 26.8 & 1.6 \\
70\_5\_10 & 3.4 & 1.2 & 3.7 & 1.9 & 4.6 & 2.5 & 3.9 & 5.0 & 0.5 \\
70\_5\_20 & 2.9 & 1.5 & 3.4 & 2.5 & 4.0 & 3.3 & 2.4 & 10.0 & 0.5 \\
70\_5\_40 & 3.4 & 1.7 & 3.5 & 2.9 & 3.6 & 4.6 & 2.4 & 11.4 & 0.5 \\
70\_5\_80 & 3.7 & 2.0 & 3.5 & 3.7 & 3.5 & 6.0 & 2.3 & 17.7 & 0.7 \\
70\_10\_10 & 2.8 & 1.3 & 3.1 & 2.1 & 3.6 & 2.6 & 2.7 & 5.0 & 0.8 \\
70\_10\_20 & 2.5 & 1.8 & 2.7 & 2.8 & 3.1 & 4.0 & 2.3 & 9.2 & 0.7 \\
70\_10\_40 & 3.3 & 2.1 & 3.2 & 3.5 & 3.5 & 5.2 & 2.3 & 10.0 & 0.8 \\
70\_10\_80 & 3.5 & 2.3 & 3.5 & 4.1 & 3.1 & 6.9 & 2.1 & 18.9 & 0.9 \\ \hline
\textbf{Average} & \textbf{3.4} & \textbf{5.2} & \textbf{3.4} & \textbf{8.1} & \textbf{3.6} & \textbf{10.8} & \textbf{2.9} & \textbf{19.1} & \textbf{1.2}
\end{tabular}%
}
\end{table}
To measure the impact of the local search procedure, we compare different strategies that differ in the frequency of execution of the local search procedure. We test three other variants of the algorithm in which the local search is called every 1, 2, and 4 iterations. The results are summarized in Table \ref{tab:LSvars} where we can observe the increased computational complexity that the local search can add to each iteration. Performing the local search procedure very frequently can lead to good solutions in fewer iterations, but it also results in longer computational times. The proposed strategy performs the local search at iterations where the reconstructed solution is better than the incumbent one, and we show that this strategy performs the best. This method allows us to perform the local search in a fewer number of iterations, but at the same time, has the potential to result in promising and better solutions.

\begin{table}[]
\centering
\caption{Performance summary of the four removal operators. The instances are grouped per number of ships.}
\label{tab:desOper}
\resizebox{\textwidth}{!}{%
\begin{tabular}{c|rrr|rrr|rrr|rrr}
\multirow{2}{*}{\textbf{\begin{tabular}[c]{@{}c@{}}Number \\ of ships\end{tabular}}} & \multicolumn{3}{c|}{\textbf{Cost-berth removal}} & \multicolumn{3}{c|}{\textbf{Cost-time removal}} & \multicolumn{3}{c|}{\textbf{Shaw removal}} & \multicolumn{3}{c}{\textbf{Random removal}} \\ \cline{2-13} 
 & \multicolumn{1}{c}{\% its} & \multicolumn{1}{c}{\begin{tabular}[c]{@{}c@{}}\% new \\ best\end{tabular}} & \multicolumn{1}{c|}{\begin{tabular}[c]{@{}c@{}}\% better \\ than current\end{tabular}} & \multicolumn{1}{c}{\% its} & \multicolumn{1}{c}{\begin{tabular}[c]{@{}c@{}}\% new \\ best\end{tabular}} & \multicolumn{1}{c|}{\begin{tabular}[c]{@{}c@{}}\% better \\ than current\end{tabular}} & \multicolumn{1}{c}{\% its} & \multicolumn{1}{c}{\begin{tabular}[c]{@{}c@{}}\% new \\ best\end{tabular}} & \multicolumn{1}{c|}{\begin{tabular}[c]{@{}c@{}}\% better \\ than current\end{tabular}} & \multicolumn{1}{c}{\% its} & \multicolumn{1}{c}{\begin{tabular}[c]{@{}c@{}}\% new \\ best\end{tabular}} & \multicolumn{1}{c}{\begin{tabular}[c]{@{}c@{}}\% better \\ than current\end{tabular}} \\ \hline
30 & 6 & 4 & 2 & 17 & 21 & 13 & 28 & 20 & 24 & 50 & 55 & 62 \\
50 & 6 & 4 & 1 & 47 & 50 & 38 & 12 & 8 & 8 & 35 & 38 & 52 \\
70 & 6 & 4 & 2 & 69 & 81 & 84 & 6 & 2 & 1 & 19 & 12 & 13
\end{tabular}%
}
\end{table}

\begin{table}[]
\centering
\caption{Performance summary of the four insertion operators. The instances are grouped per number of ships.}
\label{tab:repOper}
\resizebox{\textwidth}{!}{%
\begin{tabular}{c|rrrr|rrrr|rrrr|rrrr}
\multirow{2}{*}{\textbf{\begin{tabular}[c]{@{}c@{}}Number\\ of ships\end{tabular}}} & \multicolumn{4}{c|}{\textbf{Efficient packing insertion}} & \multicolumn{4}{c|}{\textbf{Random greedy insertion}} & \multicolumn{4}{c|}{\textbf{Arrival greedy insertion}} & \multicolumn{4}{c}{\textbf{$\kappa$-regret insertion}} \\ \cline{2-17} 
 & \multicolumn{1}{c}{\% its} & \multicolumn{1}{c}{\begin{tabular}[c]{@{}c@{}}\% of \\ time\end{tabular}} & \multicolumn{1}{c}{\begin{tabular}[c]{@{}c@{}}\% new \\ best\end{tabular}} & \multicolumn{1}{c|}{\begin{tabular}[c]{@{}c@{}}\% better \\ than current\end{tabular}} & \multicolumn{1}{c}{\% its} & \multicolumn{1}{c}{\begin{tabular}[c]{@{}c@{}}\% of \\ time\end{tabular}} & \multicolumn{1}{c}{\begin{tabular}[c]{@{}c@{}}\% new \\ best\end{tabular}} & \multicolumn{1}{c|}{\begin{tabular}[c]{@{}c@{}}\% better \\ than current\end{tabular}} & \multicolumn{1}{c}{\% its} & \multicolumn{1}{c}{\begin{tabular}[c]{@{}c@{}}\% of \\ time\end{tabular}} & \multicolumn{1}{c}{\begin{tabular}[c]{@{}c@{}}\% new \\ best\end{tabular}} & \multicolumn{1}{c|}{\begin{tabular}[c]{@{}c@{}}\% better \\ than current\end{tabular}} & \multicolumn{1}{c}{\% its} & \multicolumn{1}{c}{\begin{tabular}[c]{@{}c@{}}\% of \\ time\end{tabular}} & \multicolumn{1}{c}{\begin{tabular}[c]{@{}c@{}}\% new \\ best\end{tabular}} & \multicolumn{1}{c}{\begin{tabular}[c]{@{}c@{}}\% better \\ than current\end{tabular}} \\ \hline
30 & 5 & 2 & 4 & 1 & 17 & 10 & 24 & 7 & 10 & 5 & 13 & 3 & 68 & 83 & 60 & 89 \\
50 & 5 & 3 & 3 & 1 & 37 & 33 & 34 & 22 & 8 & 6 & 6 & 2 & 49 & 58 & 57 & 75 \\
70 & 6 & 3 & 3 & 1 & 47 & 39 & 51 & 46 & 9 & 6 & 6 & 3 & 38 & 52 & 39 & 50
\end{tabular}%
}
\end{table}
Tables \ref{tab:desOper} and \ref{tab:repOper} summarized the performance of the different removal and insertion operators used within the ALNS method. For each removal operator, we compute three metrics: (i) the percentage of iterations in which the operator was selected \textit{\% its}, (ii) the percentage of current best solutions found using the operator \textit{\% new best}, and (iii) the percentage of times that the resulting solution was better than the current one using the operator \textit{\% better than current}.
For the insertion methods, we also display a fourth column \textit{\% of time}, which indicates the percentage of time each operator has consumed from the total time spent repairing solutions. The time spent in removal methods is significantly lower than in insertion operators; therefore, we do not compute this metric for the removal operators.
We observe that the random removal is the better-performing removal method when the number of ships is 30. However, for larger instances, the cost-time removal performs better. This suggests that when the number of port visits increases, the probability of removing port visits that are not related at all also increases, making pure random methods less efficient.
Furthermore, removing port visits that overlap in time is more effective than removing the ones that overlap in berthing space. While ships can berth at multiple positions along the quay without major delays and disruptions in their schedule, berthing earlier or later can negatively impact the sailing in the rest of the voyage legs. Therefore, their flexibility comes at a larger cost. 
We also notice that the Shaw removal becomes less useful in larger instances. This operator removes pairs of port visits at a time. These pairs can be highly unrelated between them,  worsening the effects of the overall operator. Similarly to the random removal, this inter-relatedness issue increases with the number of port visits.  

Regarding the repair methods, the $\kappa$-regret insertion method shows the best and more robust performance. Even if, in some cases, it is not the most frequently used operator, it is clearly computationally intensive, and more than half of the time is spent computing the insertions related to this method. The randomized greedy insertion is also relevant and becomes more effective in larger instances. The remaining two insertion operators, packing and arrival greedy, show a modest performance. Still, they also proved to help achieve some of the best-found solutions during the algorithm run.

\begin{figure}
    \centering
    \includegraphics[width=1\textwidth]{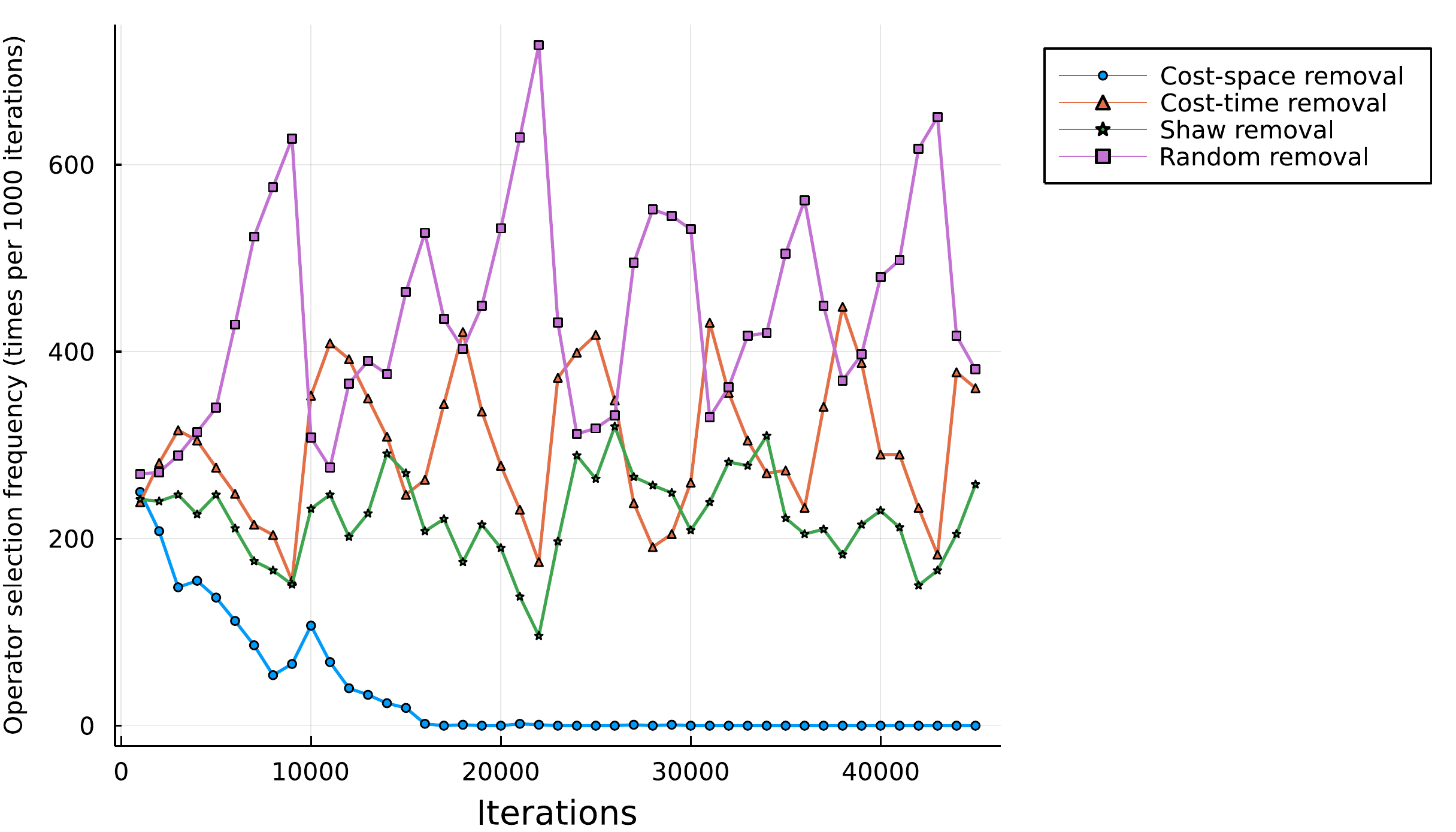}
    \caption{Usage of each removal operator during an algorithm run of 1 hour for an instance from group $30\_10\_10$.}
    \label{fig:DesOpRun}
\end{figure}
\begin{figure}
    \centering
    \includegraphics[width=1\textwidth]{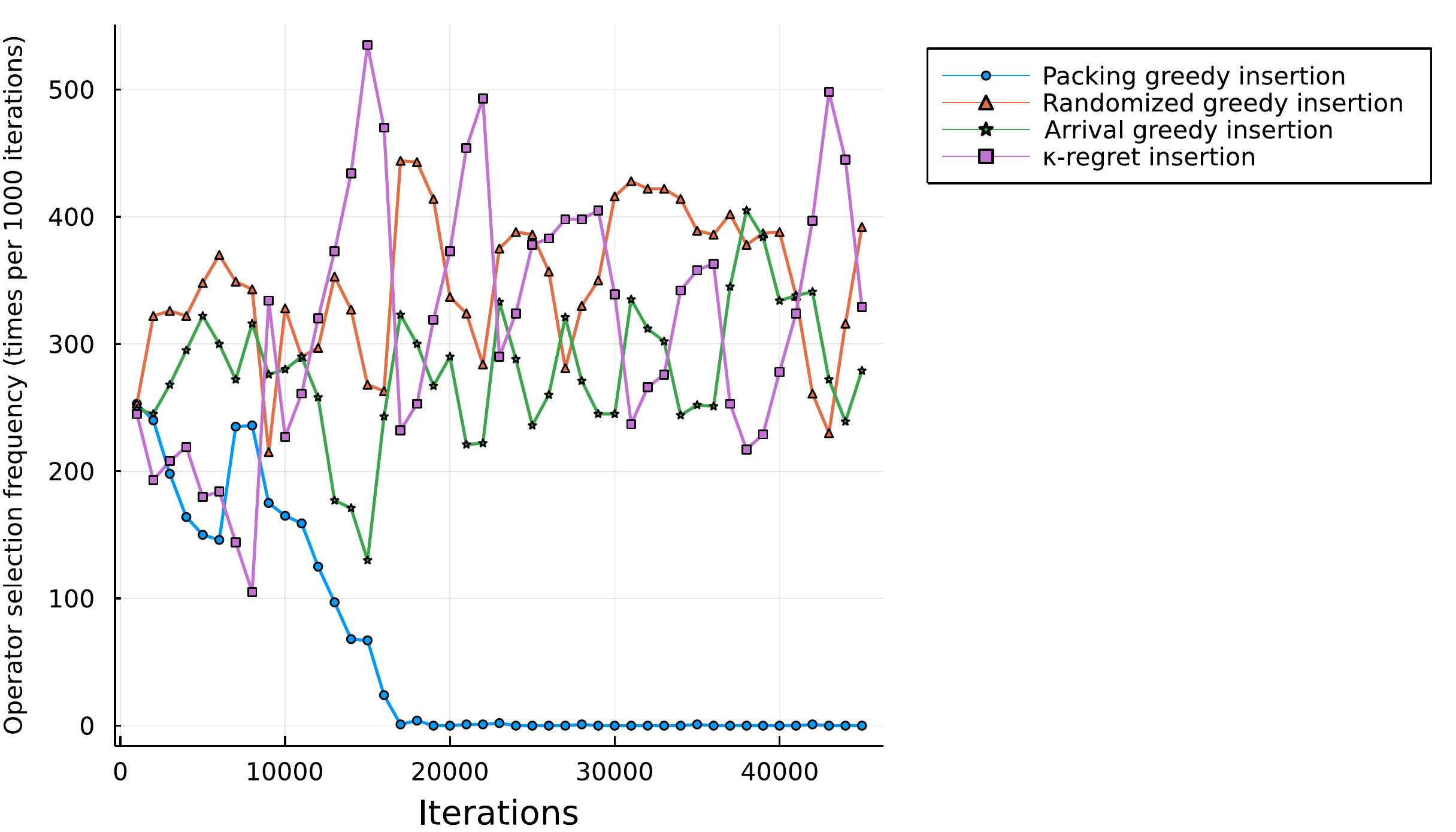}
    \caption{Usage of each insertion operator during an algorithm run of 1 hour for an instance from group $30\_10\_10$.}
    \label{fig:RepOpRun}
\end{figure}
To better understand the algorithm's behavior and when the operators are used, we tracked the use of each operator during the algorithm run. Figures \ref{fig:DesOpRun} and \ref{fig:RepOpRun} show an example run of one hour for an instance with 30 ships, 10 external ships per port, and a quay segment length of 10 meters. We observe that the cost-berth removal is only used at the beginning of the run when each operator has a more balanced probability of being chosen. The packing greedy heuristic shows a similar behavior and correlates with the low usage of these two operators as shown in Table \ref{tab:desOper} and \ref{tab:repOper}.
Nonetheless, the rest of the operators are used for most algorithm runs. Some operators show an oscillating behavior, such as the cost-time removal operator. This behavior correlates with the temperature of the acceptance criterion, which is reheated periodically and suggests that the operator has a higher probability of being selected when the temperature is higher. 

Additionally, for both removal and insertion operators, we tested the algorithm removing the worst performing operators, one at a time, but the performance of the method worsened in all cases, indicating that all operators are to some extent useful and combine well together.

\subsection{Practical impact}

This problem involves two main stakeholders, namely, the terminal operators and the shipping carriers.
The objective function covers the operational costs of both of them. This section disaggregates and analyzes the different operational costs by performing various sensitivity analyses. 

\begin{table}[]
\centering
\caption{Average operational costs and cost variation across instances with different quay segment lengths. All instances are run for one hour. The costs are in thousands of USD.}
\label{tab:impactQuaySegment}
\resizebox{0.9\textwidth}{!}{%
\begin{tabular}{r|rrrrrrr}
\multicolumn{1}{c|}{\textbf{Segment length}} & \multicolumn{1}{c}{\textbf{Waiting}} & \multicolumn{1}{c}{\textbf{Handling}} & \multicolumn{1}{c}{\textbf{Delay}} & \multicolumn{1}{c}{\textbf{Fuel}} & \multicolumn{1}{c}{\textbf{Total}} & \multicolumn{1}{c}{\textbf{Penalty}} & \multicolumn{1}{c}{\textbf{Iterations x1000}} \\ \hline
10 & 390 & 537 & 3210 & 1274 & 5411 & 0.04 & 8 \\
20 & 368 & 537 & 3246 & 1275 & 5426 & 0.05 & 15 \\
40 & 242 & 541 & 3449 & 1276 & 5509 & 0.09 & 18 \\
80 & 241 & 542 & 3910 & 1275 & 5968 & 0.14 & 35
\end{tabular}%
}
\end{table}
In Table \ref{tab:impactQuaySegment}, we group the instances per quay segment length. We observe a natural trade-off here. A shorter segment length allows a more granular set of berthing positions and, therefore, a potentially better solution quality. However, this increases the complexity of the problem, and in the case of our method, it results in fewer iterations per hour. Despite performing less than a quarter of the iterations of the instances with 80 meters segments, the method finds better solutions for the shorter-segment instances. The improvement in objective value mainly translates into shorter delays and increased waiting time. The handling and fuel costs remain similar. The vessel time windows or port calls are already pre-planned, considering a low sailing speed. This, together with the fact that fuel costs account for a large part of the total costs, results in that ships already sailing at the slowest speed in most of the solutions (see Table \ref{tab:fuelCost}). 

\begin{table}[]
\centering
\caption{Operational costs for instances grouped by different amounts of external ships per port.}
\label{tab:externalShips}
\resizebox{0.7\textwidth}{!}{%
\begin{tabular}{c|cccccc}
\textbf{\begin{tabular}[c]{@{}c@{}}\small External ships\\ \small per port\end{tabular}} & \textbf{Waiting} & \textbf{Handling} & \textbf{Delay} & \textbf{Fuel} & \textbf{Total} & \textbf{Penalty} \\ \hline
0 & 451 & 538 & 2365 & 1278 & 4632 & 0.06 \\
5 & 287 & 538 & 3022 & 1275 & 5122 & 0.04 \\
10 & 301 & 538 & 3227 & 1275 & 5341 & 0.05 \\
20 & 461 & 538 & 2718 & 1278 & 8300 & 0.11
\end{tabular}%
}
\end{table}
Another operational aspect we inspect is the impact of the external ships in the planning process. We solve the problem for instances with a different number of external ships per port, from none to twenty ships per port. The results are summarized in Table \ref{tab:externalShips}.
The results support the rationale that an increased number of external ships per port results in a more congested berth allocation and, as a result, higher operational costs.
In this case, the port congestion is reflected in the \textit{Penalty} column, which indicates the average number of port visits per instance exceeding the latest finish time. Since this type of delay is heavily penalized, improvements in this aspect significantly impact the objective value. 
These results also indicate that the level of impact of this type of collaborative problem can increase significantly when more ships are involved. 
When more ships collaborate, their potential joint savings increase and the terminal has more planning flexibility. 

As mentioned previously, fuel consumption is the main cost driver for carriers. Fuel prices have fluctuated significantly in the last two years due to the global socio-economical and political situation. We consider that ships use a very low sulfur fuel oil (VLSFO) with an estimated price of 500 USD per metric ton. However, the prices of this fuel have ranged between 200 and 1100 USD per metric ton in the last two years. Therefore, we have also tested our method using a fuel price of 200 and 1100 USD per metric ton \citep{FuelPrice}

\begin{table}[]
\centering
\caption{Average fuel consumption per ship and sailing speed  based on different fuel prices}
\label{tab:fuelCost}
\resizebox{\columnwidth}{!}{%
\begin{tabular}{r|rrrrrr}
\multicolumn{1}{c|}{\multirow{3}{*}{\textbf{\begin{tabular}[c]{@{}c@{}}Number \\ of ships\end{tabular}}}} & \multicolumn{6}{c}{\textbf{Fuel price (USD/metric tonne)}} \\ \cline{2-7} 
\multicolumn{1}{c|}{} & \multicolumn{2}{c|}{\textbf{200}} & \multicolumn{2}{c|}{\textbf{500}} & \multicolumn{2}{c}{\textbf{1100}} \\ \cline{2-7} 
\multicolumn{1}{c|}{} & \multicolumn{1}{c|}{\begin{tabular}[c]{@{}c@{}}\small average fuel consumption\\ \small(metric tonne per ship)\end{tabular}} & \multicolumn{1}{c|}{\begin{tabular}[c]{@{}c@{}}\small average \\ \small speed (knots)\end{tabular}} & \multicolumn{1}{c|}{\begin{tabular}[c]{@{}c@{}} \small average fuel consumption\\ \small(metric tonne per ship)\end{tabular}} & \multicolumn{1}{c|}{\begin{tabular}[c]{@{}c@{}}\small average \\ \small speed (knots)\end{tabular}} & \multicolumn{1}{c|}{\begin{tabular}[c]{@{}c@{}}\small average fuel consumption\\ \small(metric tonne per ship)\end{tabular}} & \multicolumn{1}{c}{\begin{tabular}[c]{@{}c@{}}\small average \\ \small speed (knots)\end{tabular}} \\ \hline
30 & \multicolumn{1}{r|}{51.04} & \multicolumn{1}{r|}{17.12} & \multicolumn{1}{r|}{50.53} & \multicolumn{1}{r|}{17.04} & \multicolumn{1}{r|}{50.37} & 17.01 \\
50 & \multicolumn{1}{r|}{51.50} & \multicolumn{1}{r|}{17.10} & \multicolumn{1}{r|}{51.11} & \multicolumn{1}{r|}{17.04} & \multicolumn{1}{r|}{50.99} & 17.02 \\
70 & \multicolumn{1}{r|}{52.14} & \multicolumn{1}{r|}{17.19} & \multicolumn{1}{r|}{51.35} & \multicolumn{1}{r|}{17.08} & \multicolumn{1}{r|}{51.07} & 17.03 \\ \hline
Average & \multicolumn{1}{r|}{51.56} & \multicolumn{1}{r|}{17.14} & \multicolumn{1}{r|}{51.00} & \multicolumn{1}{r|}{17.05} & \multicolumn{1}{r|}{50.81} & 17.02
\end{tabular}
}
\end{table}
Table \ref{tab:fuelCost} shows the average fuel consumption per ship and sailing speed, grouped by instances with the same number of ships. We observe that the average consumption can increase by more than half a metric tonne when the fuel price decreases from 500 USD per tonne to 200 USD per tonne. This difference is more prominent in instances with a large number of ships, where more ships sail marginally faster to arrive earlier at the next port to get a better service. However, when the fuel price increases above 500 USD per metric tonne, the reductions in fuel consumption are relatively small. The main explanation for this is due to the low sailing speeds in general. The fuel costs already account for a large part of the operational costs, and the solutions indicate that ships sail close to the slowest speed of 17 knots in most cases. We observe a slight increase in average sailing speed when the fuel price is low, but the size of the instance does not have an impact on the speed. A similar sensitivity analysis performed by \cite{venturini2017a} indicated a similar behavior.

\section{Conclusions}\label{sec:conclusions}

In this work, we address an emerging problem in maritime collaborative logistics that integrates the operations of both shipping carriers and terminal operators. We present both a new MIP formulation for the 
 multi-port continuous berth allocation problem with speed optimization, and an ALNS algorithm to solve it. The ALNS algorithm takes advantage of a diverse set of tailored insertion and removal methods. It guides the algorithm by prioritizing the better-performing methods. The modular characteristic of the algorithm could be exploited to develop a decision support tool for terminal operators, where the operators' experience can lead to new tailored operators. Furthermore, in terms of computational performance, the heuristic method is able to find high-quality solutions to larger instances than the ones studied in the literature and outperforms commercial solvers such as CPLEX. We also study the practical impact of the problem in terms of operational costs for the carriers and terminal operators and analyze the resulting quality of the berth plans and sailing speeds. We conclude that engaging in this type of collaboration can result in overall cost reductions for the stakeholders and also benefits to the environment due to the potential lower fuel consumption.

Some aspects of this study remain as future work or research direction.
Regarding the solution method, the insertion operators are the main bottleneck in terms of computational complexity. One could explore simpler insertion operators or other heuristic variants.
Studying the scalability of the method in more detail could be relevant. There is no doubt that the heuristic method scales better than CPLEX, and results in small instances indicate that the ALNS achieves near-optimal solutions. For larger instances, the optimality gap of CPLEX increases significantly, and the lower bound becomes impractical.
Finally, incorporating practical aspects such as transhipments or disruptions management is an attractive research direction. We envision the use of frameworks such as stochastic programming to tackle this type of problem.
All in all, this type of study highlights the potential impact of collaborative logistics and the value of integration in the transportation sector.
\vskip 2pt
\textbf{Acknowledgements:} The authors thank the Danish Maritime Fund for supporting this work.

\bibliography{sample}





\bibliographystyle{elsarticle-harv}\biboptions{authoryear}

\appendix

\end{document}